\documentclass{siamltex}

\usepackage{amscd,amsmath,amssymb,verbatim,color}

\newtheorem{problem}[theorem]{Problem}
\newtheorem{remark}[theorem]{Remark}
\newtheorem{remarks}[theorem]{Remarks}

\newcommand{\linspan}{\mathop{\rm span}\nolimits}

\newcommand{\esssup}{\mathop{\rm ess\ sup}}

\newcommand{\p}{\partial}

\newcommand{\R}{{\mathbb R}}

\newcommand{\N}{{\mathbb N}}

\newcommand{\EE}{{\mathcal E}}

\newcommand{\LL}{{\mathcal L}}

\newcommand{\WW}{{\mathcal W}}
\newcommand{\XX}{{\mathcal X}}
\newcommand{\YY}{{\mathcal Y}}

\title{Internal exponential stabilization to a non-stationary solution for 3D Navier--Stokes equations}
\author{Viorel Barbu\thanks{Department of Mathematics, University
``Al. I. Cuza'', 6600 Iasi, Romania; e-mail: {\tt vb41@uaic.ro}} 
\and S\'ergio S.~Rodrigues\thanks{University of Cergy-Pontoise, Department of Mathematics, UMR CNRS 8088, F-95000 Cergy-Pontoise; e-mail: 
{\tt Sergio.Rodrigues@u-cergy.fr}} 
\and Armen Shirikyan\thanks{University of Cergy-Pontoise, Department of Mathematics, UMR CNRS 8088, F-95000 Cergy-Pontoise; e-mail: 
{\tt Armen.Shirikyan@u-cergy.fr}}}

\begin{document}
\maketitle

\begin{abstract}
We consider the Navier--Stokes system in a bounded domain with a
smooth boundary. Given a sufficiently regular time-dependent global solution, we construct a finite-dimensional feedback control that is supported
by a given open set and stabilizes the linearized equation. The
proof of this fact is based on a truncated observability
inequality, the regularizing property for the linearized equation,
and some standard techniques of the optimal control theory. We
then show that the control constructed for the linear problem
stabilizes locally also the full Navier--Stokes system.
\end{abstract}




\pagestyle{myheadings}
\thispagestyle{plain}
\markboth{V.BARBU, S.S.RODRIGUES, AND A.SHIRIKYAN}{EXPONENTIAL STABILIZATION  FOR 3D NAVIER--STOKES EQUATIONS}

\section{Introduction}
Let $\Omega\subset\R^3$ be a connected bounded domain located locally on one side of its smooth boundary $\Gamma=\p\Omega$. We consider the
controlled Navier--Stokes system in~$\Omega$:
\begin{align}
 \p_t u+\langle u\cdot\nabla\rangle u-\nu\Delta u+\nabla p&=h+\zeta, \quad \nabla \cdot u=0,\label{01}\\
 u\bigr|_\Gamma &=0. \label{02}
\end{align}
Here $u=(u_1,u_2,u_3)$ and~$p$ are unknown velocity field and
pressure of the fluid, $\nu>0$ is the viscosity, $\langle
u\cdot\nabla\rangle$ stands for the differential operator
$u_1\p_1+u_2\p_2+u_3\p_3$, $h$~is a fixed function, and~$\zeta$ is
a control taking values in the space~$\EE$ of square-integrable
functions in~$\Omega$ whose support in~$x$ is contained in a given
open subset $\omega\subset\Omega$. The problem of exact
controllability for~\eqref{01}, \eqref{02} was in the focus of
attention of many researchers starting from the early nineties,
and it is now rather well understood. Namely, it was proved that,
given a time~$T>0$ and a smooth solution~$\hat u$ of~\eqref{01},
\eqref{02} with $\zeta\equiv0$, for any initial function~$u_0$ sufficiently close to~$\hat u(0)$ one can find a control
$\zeta:[0,T]\to\EE$ such that the solution of problem~\eqref{01},
\eqref{02} supplemented with the initial condition
\begin{equation} \label{03}
u(0,x)=u_0(x)
\end{equation}
is defined on~$[0,T]$ and satisfies the relation $u(T)=\hat u(T)$.
We refer the reader
to~\cite{cor_fur,FI1996,FE-1996,FE1999,iman01,FGIP-2004} for the exact
statements and the proofs of these results.

Even though the property of exact controllability is quite
satisfactory from the mathematical point of view, many problems
arising in applications require that the control in question be
feedback, because closed-loop controls are usually more stable
under perturbations (e.g., see the introduction to Part~3
in~\cite{coron2007}). This question has found a positive answer in the context of stabilization theory. 
It was intensively studied for the
case in which the target solution~$\hat u$ is stationary (in
particular, the external force is independent of time). A typical
result in such a situation claims that, given a smooth stationary
state~$\hat u$ of the Navier--Stokes system and a
constant~$\lambda>0$, one can construct a continuous linear
operator $K_{\hat u}:L^2\to\EE$ with finite-dimensional range such that the solution of
problem~\eqref{01}~--~\eqref{03} with $\zeta=K_{\hat u}(u-\hat u)$
and a function~$u_0$ sufficiently close to~$\hat u$ is defined for
all~$t\ge0$ and converges to~$\hat u$ at least with the
rate~$e^{-\lambda t}$. We refer the reader to the
papers~\cite{fursikov-2001,fursikov-2004,raymond-2006,raymond-2008,RT-2009,BT-2009}
for boundary stabilization and to~\cite{barbu-2003,BT-2004,BLT-2006} for stabilization by a
distributed control.

\smallskip
The aim of this paper is to establish a similar result in the case
when the target solution~$\hat u$ depends on time. Namely, we
will prove the following theorem, whose exact formulation is
given in Section~\ref{S:nonlinear}.

\medskip
{\sc Main Theorem.}
{\it Let~$(\hat u,\hat p)$ be a global smooth solution for
problem~\eqref{01}, \eqref{02} with $\zeta\equiv0$ such that
\begin{equation*}
\esssup_{(t,x)\in Q}\,\bigl|\p_t^j\p_x^\alpha \hat u(t,x)\bigr|\le
R \quad\mbox{for $j=0,1$, $|\alpha|\le 1$},
\end{equation*}
where $Q=\R_+\times\Omega$ and $R>0$ is a constant. Then for any
$\lambda>0$ and any open subset $\omega\subset\Omega$ there is an
integer $M=M(R,\lambda,\omega)\ge1$, an $M$-dimensional space
$\EE\subset C_0^\infty(\omega,\R^3)$, and a family of continuous
linear operators $K_{\hat u}(t):L^2(\Omega,\R^3)\to \EE$, $t\ge0$,
such that the following assertions hold.
\begin{itemize}
\item[\bf(a)] The function $t\mapsto K_{\hat u}(t)$ is continuous
in the weak operator topology, and its operator norm is bounded
by a constant depending only on~$R$, $\lambda$, and~$\omega$.
\item[\bf(b)] For any divergence free function $u_0\in
H_0^1(\Omega,\R^3)$ that is sufficiently close to~$\hat u(0)$ in
the $H^1$-norm problem~\eqref{01}~--~\eqref{03} with~$\zeta=K_{\hat u}(t)(u-\hat u(t))$ has a unique global strong
solution $(u,p)$, which satisfies the inequality
\begin{equation*}
|u(t)-\hat u(t)|_{H^1}\le Ce^{-\lambda t}|u_0-\hat
u(0)|_{H^1}, \quad t\ge0.
\end{equation*}
\end{itemize}
}

Note that this theorem remains true for the two-dimensional
Navier--Stokes system, and in this case, it suffices to assume
that the initial function~$u_0$ is close to~$\hat u(0)$ in the
$L^2$-norm. Furthermore, the approach developed in this
paper applies equally well to the case when the control acts via the
boundary. This situation will be addressed in a subsequent publication.

\smallskip
As was mentioned above, the problem of feedback stabilization is rather well understood for stationary reference solutions. Let us explain informally the additional difficulties arising in the non-stationary case and reveal a common mechanism of stabilization. 
A well-known argument based on the contraction mapping principle
enables one to prove that a control stabilizing the linearized
problem locally stabilizes also the nonlinear equation. Thus, it suffices to study the linearized problem. The main idea in the stationary case is to split it into a system of two autonomous equations, the first of which is finite-dimensional and has the zero solution as a possibly unstable equilibrium point, whereas the second is exponentially stable due to the large negative eigenvalues of the Laplacian. One then applies methods of finite-dimensional theory (e.g., the pole assignment theorem~\cite[Section~2.5]{zabczyk1992})   to find a stabilizing feedback control for the first equation and proves that using the same control in the original problem yields an exponentially decaying solution; see~\cite{barbu-2003,BT-2004,raymond-2006,BLT-2006,raymond-2008,RT-2009,BT-2009}.

It is difficult to apply this approach in the case of time-dependent reference solutions, because a non-autonomous equation does not necessarily admits invariant subspaces. However, the above-mentioned scheme for stabilization is based essentially on the so-called Foia\c s--Prodi property for parabolic PDE's~\cite{FP-1967}. It says, roughly speaking, that if the projections of  two solutions to the unstable modes   converge to each other as time goes to infinity, then the difference between these solutions goes to zero. It turns out that the conclusion remains true if the projections are close to each other at times proportional to a fixed constant. The main  idea of this paper is to choose a control that ensures the equality at integer times for the projections of two solutions to the unstable modes. More precisely, we consider the following problem obtained
by linearizing~\eqref{01}, \eqref{02} around a non-stationary solution~$\hat u(t,x)$:
\begin{equation} \label{06}
 \p_t v+\langle \hat u \cdot\nabla\rangle v+\langle v\cdot\nabla\rangle\hat u -\nu\Delta v+\nabla p=\zeta, \quad \nabla \cdot v=0, \quad v\bigr|_\Gamma =0.
\end{equation}
Let us assume that, for a sufficiently large integer~$N$, we have
constructed a continuous linear operator $\bar
\zeta:L^2(\Omega,\R^3)\to L^2((0,1);\EE)$ such that, for any
initial function~$v_0$, the solution of~\eqref{06} with
$\zeta=\bar\zeta(v_0)$ issued from~$v_0$ satisfies the relation
$\Pi_N v(1)=0$, where~$\Pi_N$ stands for the orthogonal projection
in~$L^2$ onto the subspace spanned by the first~$N$ eigenfunctions of
the Stokes operator in~$\Omega$. In this case, using the
Poincar\'e inequality and regularizing property of the resolving
operator for~\eqref{06}, we get
\begin{align}
|v(1)|_{L^2}&=|(I-\Pi_N)v(1)|_{L^2}\le C_1\alpha_N^{-1/2}|v(1)|_{H^1}\notag\\
&\le
C_2\,\alpha_N^{-1/2}\bigl(|v_0|_{L^2}+|\bar\zeta(v_0)|_{L^2((0,1);\EE)}\bigr)
\le C_3\,\alpha_N^{-1/2}|v_0|_{L^2},\label{07}
\end{align}
where $\{\alpha_j\}$ denotes the increasing sequence of the
eigenvalues for the Stokes operator and $C_i$, $i=1,2,3$, are some
constants not depending on~$N$. The fact that~$C_3$ is
independent of~$N$ is a crucial property, and its proof is based
on a truncated observability inequality (see Section~\ref{s5.3} in
Appendix). It follows from~\eqref{07} that, if~$N$ is sufficiently
large, then $|v(1)|_{L^2}\le e^{-\lambda}|v_0|_{L^2}$.
Iterating this procedure, we get an exponentially decaying
solution. Once an exponential stabilization of the linearized
problem~\eqref{06} is obtained, the existence of an exponentially stabilizing feedback
control can be proved with the help of the dynamic programming
principle. We refer the reader to Section~\ref{S:Lin} for an
accurate presentation of the results on the linearized equation and some further comments on the existence of Lyapunov function, derivation of a Riccati equation for the feedback control operator, and the dimension of controllers.

\smallskip
The paper is organized as follows. In Section~\ref{S:prelim}, we
introduce the functional spaces arising in the theory of the
Navier--Stokes equations and recall some well-known facts.
Section~\ref{S:Lin} is devoted to studying the linearized problem.
In Section~\ref{S:nonlinear}, we establish the main result of the
paper on local exponential stabilization of the full
Navier--Stokes system. The Appendix gathers some auxiliary results
used in the main text.

\subsection*{Notation}
We write~$\mathbb N$ and~$\mathbb R$ for the sets of non-negative integers and real numbers,
respectively, and we define $\mathbb N_0=\mathbb N\setminus\{0\}$, $\R_+=(0,\,+\infty)$, and $\R_s=(s,+\infty)$. We denote by $\Omega\subset\R^3$ a bounded domain with a $C^2$-smooth boundary $\Gamma=\p\Omega$, and for $\tau\in\R$, we set $I_\tau=(\tau,\tau+1)$,  $Q_\tau=I_\tau\times\Omega$, and $Q=\R_+\times \Omega$. 
The partial time derivative~$\partial_tu$ of a function $u(t,x)$ will be denoted by~$u_t$.

\smallskip
For a Banach space~$X$, we denote by~$|\cdot|_X$ the corresponding norm, by~$X'$ its dual, and by~$\langle \cdot,\cdot\rangle_{X',X}$ the duality between~$X'$ and~$X$.

\smallskip
If  $X$ and~$Y$ are Banach spaces and $I\subseteq\mathbb R$ is an open interval, then we write
$$
W(I,\,X,\,Y):=\{f\in L^2(I,\,X)\mid\,f_t\in L^2(I,\,Y)\},
$$
where the derivative $f_t=\frac{df}{dt}$ is taken in the sense of distributions. 
This space is endowed with the natural norm
$$
|f|_{W(I,\,X,\,Y)}:=\bigl(|f|_{L^2(I,\,X)}^2+|f_t|_{L^2(I,\,Y)}^2\bigr)^{1/2}.
$$
Note that if $X=Y$, then we obtain the Sobolev space $W^{1,2}(I,X)$.

\smallskip
If $I\subset\R$ is a closed interval, then~$C(I,X)$ stands for the space of continuous functions $f:I\to X$ with the norm
$$
|f|_{C(I,X)}=\max_{t\in I}|f(t)|_X. 
$$

\smallskip
For a given space~$Z$ of functions $f=f(t)$ defined on an interval of~$\mathbb R$ and a constant~$\lambda>0$, we define
$$
Z_\lambda:=\{f\in Z\mid\,e^{(\lambda/2) t}f\in Z\}.
$$
This space is endowed with the norm
$$
|f|_{Z_\lambda}:=\bigl(|f|_Z^2+|e^{(\lambda/2) t}f|_Z^2\bigr)^{1/2}.
$$

\smallskip
Throughout the paper, we deal with two integers, $N$ and~$M$. Roughly speaking, $N$ stands for the number of unstable modes in the linearized Navier--Stokes system and~$M$ denotes the space dimension of the control function arising in various problems. 

\smallskip
$\overline C_{[a_1,\dots,a_k]}$ denotes a function of non-negative variables~$a_j$ that increases in each of its arguments.

\smallskip
$C_i$, $i=1,2,\dots$, stand for unessential positive constants.

\section{Preliminaries}\label{S:prelim}

\subsection{Functional spaces and reduction to an evolution equation}
In what follows, we will confine ourselves to the 3D case, although all the results remain valid for the 2D Navier--Stokes equations.

Let $\Omega\subset \mathbb R^3$ be a connected bounded domain located
locally on one side of its $C^2$-smooth boundary $\Gamma=\partial \Omega$.
It is natural to study the incompressible Navier--Stokes system as
an evolution  equation in the subspace~$H$ of divergence free
vector fields tangent to the boundary:
$$
H:=\{u\in L^2(\Omega,\R^3)\mid\,\nabla\cdot u=0\text{ in }\Omega,\,
u\cdot\mathbf n=0 \text{ on }\Gamma\}.
$$
Here $L^2(\Omega,\R^3)$ is the space of square
integrable vector fields $(u_1,u_2,u_3)$ in~$\Omega$, $\nabla\cdot
u:=\partial_1u_1+\partial_2u_2+\partial_3u_3$ is the  divergence
of~$u$, and~$\mathbf n$ is the normal vector to
the boundary~$\Gamma$. Let us denote by $H^s(\Omega)$ the Sobolev space of order~$s$ and by~$H^s(\Omega,\R^3)$ the space of vector fields in~$\Omega$ whose components belong to~$H^s(\Omega)$. To simplify notation, we will often write~$L^2$ and~$H^s$; the context will imply the domain on which these spaces are considered. Define 
$$
V:=\{u\in H^1(\Omega,\R^3)\mid\,\nabla\cdot u=0\text{ in }\Omega,\, u=0 \text{ on }\Gamma\}, \quad
U:=H^2(\Omega,\R^3)\cap V.
$$
Note that~$U$ coincides with the natural domain~${\rm D}(L)$ of the Stokes operator $L=-\nu\Pi\Delta$,
where~$\Pi$ is the orthogonal projection in~$L^2(\Omega,\R^3)$ onto~$H$. The spaces~$H$, $V$ and~$U$ are endowed with the scalar products
$$
(u,\,v)_H:=(u,\,v)_{L^2(\Omega,\R^3)},\quad
(u,\,v)_V:=\langle Lu,\,v\rangle_{V^\prime,V},\quad
(u,\,v)_{{\rm D}(L)}:=(Lu,\,Lv)_{L^2(\Omega,\R^3)},
$$
respectively, and we denote by~$|\cdot|_H$, $|\cdot|_V$ and
$|\cdot|_{{\rm D}(L)}$ the corresponding norms. Finally, for any integer $k\ge0$, we introduce the space Banach~$\mathcal W^k$ of measurable vector functions $u=(u_1,u_2,u_3)$ defined in~$Q$ such that 
\begin{equation} \label{Wspace}
|u|_{\WW^k}:=\sum_{j,\alpha}\,
\esssup_{(t,x)\in Q}\,\bigl|\p_t^j\p_x^\alpha u(t,x)\bigr|<\infty,
\end{equation}
where the sum is taken over $0\le j\le k$ and $|\alpha|\le 1$. In the case $k=1$, we will write~$\WW$ instead of~$\WW^1$. 

\medskip
It is well known (e.g., see~\cite{temams}) that problem~\eqref{01}, \eqref{02} is equivalent to the following evolutionary equation in~$H$:
\begin{equation} \label{00}
u_t + Lu + B u=\Pi(h+\zeta),
\end{equation}
where $Bu:V\to V^\prime$ is defined by $Bu:=B(u,\,u)$ with
$$
\langle B(u,\,v),\,w\rangle_{V^\prime,V}
=\sum_{i,j=1}^3\int_\Omega u_i(\partial_i u_j)w_j{\rm d}x.
$$
In the following, we will deal also with linear equations obtained from~\eqref{00}
after replacing~$B$ by one of the operators~$\mathbb B(\hat u)$ and~$\mathbb
B^\ast(\hat u)$, where $\hat u\in\WW$ is a fixed function, $\mathbb B(\hat u)v=B(v,\hat u)+B(\hat u,v)$, and~$\mathbb B^\ast(\hat u)$ stands for the formal adjoint of~$\mathbb B(\hat u)$ with respect to the scalar product on~$H$:
$$
\langle \mathbb B^*(\hat u)v,\,w\rangle_{V^\prime,V}
=\sum_{i,j=1}^3\int_\Omega (v_j\p_i\hat u_j-\hat u_j\p_j v_i)w_i{\rm d}x.
$$
Namely, let us consider the problem
\begin{align}
r_t+Lr+\hat{\mathbb B}r&=f,\quad t\in I_0=(0,\,1),\label{eq.r}\\
r(0)&=r_0,\label{ic.r}
\end{align}
where $\hat{\mathbb B}=\mathbb B(\hat u)$ or $\mathbb B^*(\hat u)$.

\begin{lemma}\label{bd.hatu-r0}
For any $\hat u\in\WW$, $u_0\in H$, and $f\in L^2(I_0,\,V')$, problem~\eqref{eq.r}, \eqref{ic.r} has a unique solution $r\in W(I_0,V,V')$,
which satisfies the inequality
\begin{equation} \label{l1}
|r|_{C(\bar I_0,\,H)}^2+\int_{I_0}|r|_V^2\,{\rm d}t+\int_{I_0}|r_t|_{V^\prime}^2\,{\rm d}t
\leq\overline C_{[|\hat u|_{\WW^0}]} |r_0|_H^2+|f|_{L^2(I_0,V')}^2,
\end{equation}
where $\bar I$ stands for the closure of an interval~$I\subset\R$. 
Moreover, if $f\in L^2(I_0,H)$, then we have the inclusions $\sqrt{t}v\in C(\bar I_0,\,V)$, $\sqrt{t}v\in L^2(I_0,\,U)$ and the estimate
\begin{equation} \label{l2}
|\sqrt{t}r|_{C(\bar I_0,\,V)}^2+\int_{I_0}(\sqrt{t}|r|_U)^2\,{\rm
d}t \leq\overline C_{[|\hat u|_{\WW^0}]}(|r_0|_H^2+|f|_{L^2(I_0,H)}^2).
\end{equation}
Finally, if $r_0\in V$ and $f\in L^2(I_0,H)$, then $r\in W(I_0,U,H)$ and
\begin{equation} \label{l3}
|r|_{C(\bar I_0,\,V)}^2+\int_{I_0}|r|_{{\rm D} (L)}^2\,{\rm
d}t+\int_{I_0}|r_t|_{H}^2\,{\rm d}t\leq\overline C_{[|\hat
u|_{\WW^0}]}(|r_0|_V^2+|f|_{L^2(I_0,H)}^2).
\end{equation}
\end{lemma}

The proof of this lemma is based on a well-known argument, and we
will not present it here. We refer the reader to the books~\cite{lions69,temams} for more general results on existence, uniqueness, and a priori estimates for solutions linear and nonlinear Navier--Stokes type problems. 

\subsection{Setting of the problem}
Let us fix a function $h\in L^2(\R_+,\,H)$ and suppose that
$\hat u\in L^2(\R_+,\,V)\cap \mathcal W$ solves the Navier--Stokes system
$$
\hat u_t + L\hat u + B\hat u=h,\quad t> 0.
$$
Given a function $u_0\in H$ and a sub-domain $\omega\subseteq\Omega$, our
goal is to find a finite-dimensional subspace
$\mathcal E\subset L^2(\omega,\R^3)$ and a control
$\zeta\in L_{\rm loc}^2(\R_+,\,\mathcal E)$ such that the solution of
the problem
\begin{equation} \label{CP}
u_t + L u + B u=h+\Pi\zeta,\quad u(0)=u_0
\end{equation}
is defined for all $t>0$ and converges exponentially to~$\hat u$, i.e.,
$$
|u(t)-\hat u(t)|_H\leq C\,e^{-\kappa t}\quad\mbox{for $t\ge0$},
$$
where $C$ and $\kappa$ are positive constants.

\medskip
Let us write $L^2(\Omega,\R^3)$ as  a direct sum $L^2(\Omega,\R^3)=H\oplus H^\bot$,
where~$H^\bot$ denotes the orthogonal complement of~$H$ in~$L^2$.
For each positive integer $N$, we now define $N$-dimensional spaces $E_N\subset L^2$
and $F_N\subset H$ as follows. Let $\{\phi_i\mid\,i\in\mathbb N_0\}$ be an orthonormal
basis in $L^2(\Omega,\R^3)$ formed by the eigenfunctions of the Dirichlet Laplacian
and let $0<\beta_1\leq\beta_2\leq\dots$ be the corresponding eigenvalues. Furthermore, let $\{e_i\mid\,i\in\mathbb N_0\}$ be the orthonormal basis in~$H$ formed by the eigenfunctions of the Stokes operator and let $0<\alpha_1\leq\alpha_2\leq\dots$ be the corresponding eigenvalues.
For each $N\in\mathbb N_0$, we introduce the $N$-dimensional subspaces
$$
E_N:=\linspan\{\phi_i\mid\,i\leq N\}\subset L^2(\Omega,\R^3), \quad F_N:=\linspan\{e_i\mid\,i\leq N\}\subset H
$$
and denote by $P_N:L^2(\Omega,\R^3)\to E_N$ and $\Pi_N:L^2(\Omega,\R^3)\to
F_N$ the corresponding orthogonal projections. We will show that
the required control space can be chosen in the form $\EE_M=\chi
E_M$, where $\chi\in C_0^\infty(\Omega)$ is a given function not
identically equal to zero, and the integer~$M$ is sufficiently
large. 

\medskip
Let us note that, seeking a solution of~\eqref{CP} in the form $u=\hat u+v$, we obtain the following equivalent problem for~$v$:
\begin{equation} \label{ECP}
v_t + L v + B v+\mathbb B(\hat u)v=\Pi\zeta,\quad v(0)=v_0,
\end{equation}
where $v_0=u_0-u(0)$. It is clear that it suffices to consider the problem of exponential stabilization to zero for solutions of~\eqref{ECP}. Thus, in what follows, we  will study problem~\eqref{ECP}.

\section{Main result for linearized system}\label{S:Lin}

We fix a function $\hat u\in L_{\rm loc}^2(\R_+,\,V)\cap\mathcal W$. In what follows, it will be convenient to write the control~$\zeta$ entering~\eqref{ECP} in the form $\zeta=\chi P_M\eta$, where~$\eta$ takes its values in~$L^2(\Omega,\R^3)$ and~$\chi\in C_0^\infty(\Omega)$ is a nonzero function not identically equal to zero. Thus, we study the problem
\begin{align}
v_t + Lv+\mathbb B(\hat u)v &=\Pi(\chi P_M\eta),\label{eq.v}\\
v(0)&=v_0,\label{v0}
\end{align}
where $v_0\in H$. We refer the reader to~\cite{lions69,temams} for precise definitions of the concept of a solution for~\eqref{eq.v} (and all other Navier--Stokes type PDE's).

\subsection{Existence of a stabilizing control}

We begin with the following result, which shows that one can choose a finite-dimensional control exponentially stabilizing the zero solution for~\eqref{eq.v}. 

\begin{theorem}\label{Th:main.lin}
For each $v_0\in H$ and $\lambda>0$, there is an integer
$M=\overline C_{[|\hat u|_{\mathcal W},\lambda]}\ge1$ and a control
$\eta^{\hat u,\lambda}(v_0)\in L^2(\R_+,\,E_{M})$ such that the
solution~$v$ of system~\eqref{eq.v}, \eqref{v0} satisfies the inequality
\begin{equation} \label{ed1}
|v(t)|_H^2\leq\kappa_1 |v_0|_H^2 e^{-\lambda t}, \quad t\ge0,
\end{equation}
where $\kappa_1=\overline C_{[|\hat u|_{\mathcal W},\lambda]}>0$ is a constant not depending on~$v_0$.
Moreover, the mapping $v_0\mapsto \eta^{\hat u,\lambda}(v_0)$ is linear and satisfies the inequality
\begin{equation} \label{3.4}
\bigl|e^{(\tilde\lambda/2)t}\eta^{\hat u,\lambda}(v_0)\bigr|_{L^2(\R_+,E_M)}
\le \kappa_2|v_0|_H, 
\end{equation}
for $0\leq\tilde\lambda<\lambda$, where $\kappa_2=\overline{C}_{[|\hat u|_\WW,\lambda,(\lambda-\tilde\lambda)^{-1}]}$. 
Finally, if $v_0\in V$, then
\begin{equation} \label{ed2}
|v(t)|_V^2\leq\kappa_3 |v_0|_V^2 e^{-\lambda t}, \quad t\ge0,
\end{equation}
where $\kappa_3=\overline C_{[|\hat u|_{\mathcal W},\lambda]}>0$ does not depend on~$v_0$.
\end{theorem}

To prove this theorem, we will need two auxiliary lemmas. For each $\tau\geq 0$, consider equation~\eqref{eq.v} on the
time interval $I_{\tau}=(\tau,\tau+1)$ and supplement it with the initial condition
\begin{equation}
v(\tau)=w_0. \label{vt0}
\end{equation}
Let us denote by $S_{\hat u,\tau}(w_0,\,\eta)$ the operator that takes the pair~$(w_0,\,\eta)$
to the solution of~\eqref{eq.v}, \eqref{vt0}. By Lemma~\ref{bd.hatu-r0},  the operator~$S_{\hat u,\tau}$ is continuous
from~$H$ to $C(\bar I_{\tau},\,H)\cap L^2(I_{\tau},\,V)$ and from~$V$ to $C(\bar I_{\tau},\,V)\cap L^2(I_{\tau},\,U)$.
We will write $S_{\hat u,\tau}(w_0,\,\eta)(t)$ for the value of the solution at time~$t$.

\begin{lemma}\label{L:exist}
For each $N\in\mathbb N$ there is an integer $M_1=\overline C_{[\lambda,|\hat u|_{\mathcal W},N]}\ge1$ such that, for every $w_0\in H$, one can find a control $\eta\in L^2(I_{\tau},\,E_{M_1})$ for which 
$$
\Pi_N S_{\hat u,\tau}(w_0,\,\eta)(\tau+1)=0.
$$
Moreover, there is a constant $C_\chi$ depending only on~$|\hat u|_\WW$ {\rm(}but not on~$N$ and~$\tau${\rm)} such that
\begin{equation} \label{14}
|\eta|_{L^2(Q_\tau)}^2\leq C_\chi|w_0|_H^2.
\end{equation}
\end{lemma}

\begin{proof}
Let us fix $\epsilon>0$ and consider the following minimization problem.

\begin{problem}\label{probl.ep}
Given $M,N\in\mathbb N$ and $w_0\in H$, find the minimum of the quadratic functional
$$
J_\epsilon(v,\,\eta):=|\eta|_{L^2(Q_{\tau},\,\R^3)}^2
+\frac{1}{\epsilon}|\Pi_N
S_{\hat u,\tau}(w_0,\,\eta)(\tau+1)|_H^2
$$
on the set of functions $(v,\,\eta)\in W(I_\tau,\,V,\,V^\prime)\times L^2(Q_\tau,\,\R^3)$ that satisfy~\eqref{eq.v} and~\eqref{vt0}.
\end{problem}

Theorem~\ref{T:Jmin-A} implies that Problem~\ref{probl.ep} has a unique minimizer~$(\bar v_\epsilon,\bar \eta_\epsilon)$, which linearly depends on~$w_0\in H$. We now derive some estimates for the norm of the optimal control~$\bar \eta^\epsilon$.

To this end, the general theory of linear-quadratic optimal control problems is applicable. We use here a version of the Karush--Kuhn--Tucker theorem (see Theorem~\ref{T:Jmin-F}). Let us define the affine mapping
\begin{align*}
 F:W(\bar I_{\tau},\,V,\,V^\prime)\times L^2(Q_{\tau},\R^3)&\to H\times
 L^2(I_{\tau},\,V^\prime),\\
 \bigl(v,\,\eta)&\mapsto(v(0)-w_0,\,v_t+Lv+\mathbb B(\hat u)v-\Pi(\chi P_M\eta)\bigr)
\end{align*}
and note that its derivative is surjective.
Hence, by the Karush--Kuhn--Tucker theorem, there is a Lagrange multiplier
$(\mu^\epsilon,\, q^\epsilon)\in H\times L^2(I_{\tau},\,V)$ such that\,\footnote{The space $H\times L^2(I_{\tau},\,V)$
is regarded as the dual of $H\times L^2(I_{\tau},\,V')$, and the sign~$\circ$ stands for the composition of two linear operators.}
$$
J_\epsilon^\prime(\bar
v^\epsilon,\,\bar\eta^\epsilon)
-(\mu^\epsilon,\,q^\epsilon)\circ F^\prime(\bar v^\epsilon,\,\bar\eta^\epsilon)=0.
$$
It follows that, for all $(z,\,\xi)\in W(\bar I_{\tau},\,V,\,V^\prime)\times L^2(Q_{\tau},\,\R^3)$,
we have
\begin{gather}
\frac{2}{\epsilon}(\Pi_N\bar v^\epsilon(\tau+1),\,z(\tau+1))_H
+(z(\tau),\,\mu^\epsilon)_H
+\int_{I_{\tau}}\langle z_t+Lz+{\mathbb B}(\hat u)z,\, q^\epsilon\rangle_{V^\prime,V}\,{\rm d}t=0,\label{eps.z}\\
2\int_{I_{\tau}}(\bar\eta^\epsilon,\,\xi)_{L^2}\,{\rm d}t
+\int_{I_{\tau}}\langle-\Pi(\chi P_M\xi),\,q^\epsilon\rangle_{V^\prime,V}\,{\rm d}t=0.\label{eps.xi}
\end{gather}
Relation~\eqref{eps.z} implies that $ q^\epsilon$ is the solution of the problem
\begin{align}
q^\epsilon_t-Lq^\epsilon-\mathbb B^\ast(\hat u)q^\epsilon&=0,\quad  t\in I_{\tau},\label{eq.q}\\
q^\epsilon(\tau+1)&=-2\epsilon^{-1} \Pi_N\bar v^\epsilon(\tau+1).\label{q1}
\end{align}
Furthermore, it follows from~\eqref{eps.xi} that
\begin{equation}\label{bar.eta.q}
2\bar\eta_\epsilon= P_M(\chi q^\epsilon).
\end{equation}
Combining~\eqref{eq.v}, \eqref{eq.q}, and~\eqref{bar.eta.q}, we derive
\begin{align*}
\frac{\rm d}{{\rm d}t}(q^\epsilon,\,\bar
v^\epsilon)_H&=(q^\epsilon_t,\,\bar
v^\epsilon)_H+(q^\epsilon,\,\bar
v^\epsilon_t)_H\\
&=(Lq^\epsilon+\mathbb B^\ast(\hat u)q^\epsilon,\,\bar
v^\epsilon)_H+(q^\epsilon,\,-L\bar v^\epsilon-\mathbb
B(\hat u)\bar v^\epsilon+\Pi(\chi P_M\bar\eta^\epsilon))_H\\
&=(q^\epsilon,\,\Pi(\chi P_M\bar\eta^\epsilon))_H
=\frac{1}{2}|P_M(\chi q^\epsilon)|_{L^2}^2.
\end{align*}
Integrating in time over the interval~$I_\tau$, we obtain
$$
\int_{I_{\tau}}\bigl|P_M(\chi q^\epsilon(t))\bigr|_{L^2}^2\,{\rm d}t=
2\bigl((q^\epsilon(\tau+1),\,\bar
v^\epsilon(\tau+1))_H-(q^\epsilon(\tau),\,\bar
v^\epsilon(\tau))_H\bigr).
$$
Recalling now~\eqref{q1}, we see that $2(q^\epsilon(\tau+1),\,\bar
v^\epsilon(\tau+1))_H=-\epsilon|q^\epsilon(\tau+1)|_H^2$ and therefore
\begin{equation} \label{19}
\int_{I_{\tau}}|P_M(\chi q^\epsilon)|_{L^2}^2\,{\rm d}t+
\epsilon|q^\epsilon(\tau+1)|_H^2 =-2(q^\epsilon(\tau),\,\bar
v^\epsilon(\tau))_H.
\end{equation}
We wish to use the truncated observability inequality~\eqref{obs-in.M} to estimate the right-hand side of~\eqref{19}. To this end, we take $M=M_1$, where~$M_1$ is the integer constructed in Proposition~\ref{obs-in.trunc}.
Then, for every $\alpha>0$, we can write
\begin{align*}
\int_{I_{\tau}}|P_M(\chi q^\epsilon)|_{L^2}^2\,{\rm d}t
+\epsilon|q^\epsilon(\tau+1)|_H^2&\leq\alpha|q^\epsilon(\tau)|_H^2+\alpha^{-1}|\bar
v^\epsilon(\tau)|_H^2\\
&\leq \alpha D_\chi\int_{I_{\tau}}|P_M\chi
q^\epsilon|_{L^2}^2\,{\rm d}t +\alpha^{-1}|\bar
v^\epsilon(\tau)|_H^2.
\end{align*}
Setting $\alpha=(2D_\chi)^{-1}$, we obtain
\begin{equation}\label{bds.epsilon}
\int_{I_{\tau}}|P_M(\chi q^\epsilon)|_{L^2}^2\,{\rm d}t+
2\epsilon|q^\epsilon(\tau+1)|_H^2 \leq 4D_\chi|w_0|_H^2.
\end{equation}
In particular, the family of functions $\{P_M(\chi
q^\epsilon)\mid\,\epsilon>0\}$ is bounded
in~$L^2(Q_\tau,\R^3)$, and the family of solutions
$\left\{\bar v^\epsilon\mid\,\epsilon>0\right\}$ for
problem~\eqref{eq.v}, \eqref{vt0} is bounded in
$L^2(I_{\tau},\,V)$. It follows that the family $\left\{\bar
v^\epsilon_t\mid\,\epsilon>0\right\}$ is bounded in
$L^2(I_{\tau},\,V^\prime)$. Thus, we can find a sequence
$\epsilon_n\to0^+$  such that
\begin{align*}
\eta^{\epsilon_n}=\frac12 P_M(\chi q^{\epsilon_n})\quad&\rightharpoonup\quad
\eta^0\quad\text{in}\quad L^2(I_{\tau},\,E_M),\\
\bar v^{\epsilon_n}\quad&\rightharpoonup\quad v^0\quad\text{in}\quad L^2(I_{\tau},\,V),\\
\bar v^{\epsilon_n}_t\quad&\rightharpoonup\quad v_t^0\quad\text{in} \quad L^2(I_{\tau},\,V^\prime),
\end{align*}
where $\eta^0\in L^2(I_{\tau},\,E_M)$ and $v^0\in W(I_\tau,V,V')$ are some functions.
A standard limiting argument shows that~$v^0$ is a solution of problem~\eqref{eq.v}, \eqref{vt0} with $\eta=\eta^0$.
Furthermore, it follows from~\eqref{bds.epsilon} and~\eqref{q1} that
$$
|\Pi_N\bar v^\epsilon(\tau+1)|_H^2=\frac{\epsilon^2}{4} |q^\epsilon(\tau+1)|_H^2
\le \frac{\epsilon D_\chi}{2}|w_0|_H^2\to0
\quad\mbox{as}\quad \epsilon\to0.
$$
This convergence implies that $\Pi_Nv^0(\tau+1)=0$. Furthermore, it follows from~\eqref{bds.epsilon} that the function~$\eta^0$ satisfies inequality~\eqref{14} with $C_\chi=4D_\chi$. The proof of the lemma is complete.
\end{proof}

In view of Lemma~\ref{L:exist}, it makes sense to consider the following minimization problem.

\begin{problem}\label{probl}
Given integers $M, N\ge1$ and a function $w_0\in H$, find the minimum of the quadratic functional
$$
J(\eta):=|\eta|_{L^2(Q_{\tau},\,\R^3)}^2
$$
on the set of functions $(v,\,\eta)\in W(I_\tau,\,V,\,V^\prime)\times
L^2(I_\tau,\,E_M)$ satisfying equations~\eqref{eq.v}, \eqref{vt0} and the condition
$\Pi_N v(\tau+1)=0$.
\end{problem}

The following result shows that the control~$\eta$ constructed in Lemma~\ref{L:exist} can be chosen to be a linear function of the initial state. 

\begin{lemma}\label{L:min-lin.cont}
Let $N\ge1$ be an arbitrary integer and let $M$ be the  integer constructed  in Lemma~\ref{L:exist}. Then for any $w_0\in H$ Problem~\ref{probl}
has a unique minimizer $(\bar v^{\hat u,\tau},\,\bar\eta^{\hat u,\tau})\in W(I_\tau,\,V,\,V^\prime)\times L^2(I_\tau,\,E_M)$.
Moreover, the mapping $w_0\mapsto (\bar v^{\hat u,\tau},\,\bar\eta^{\hat u,\tau})$ is linear and continuous in the corresponding spaces,
and there is a constant~$C_\chi$ depending only on~$|\hat u|_\WW$ {\rm(}but not on~$N$ and~$\tau${\rm)} such that
\begin{equation} \label{22}
|\bar\eta^{\hat u,\tau}|_{L^2(I_{\tau},\,E_M)}^2\leq C_\chi|w_0|_H^2.
\end{equation}
\end{lemma}

\begin{proof}
Let us fix $N\in\mathbb N$, set
$$
W_N(I_{\tau},\,V,\,V^\prime):=\{v\in W(I_{\tau},\,V,\,V^\prime)\mid\,\Pi_Nv(\tau+1)=0\}, 
$$
and define~$\mathcal X$ as the space of functions $(v,\eta)\in W_N(I_{\tau},\,V,\,V^\prime)\times L^2(I_\tau,E_M)$ that satisfy equation~\eqref{eq.v}. In view of Lemma~\ref{L:exist} and the linearity of~\eqref{eq.v}, $\XX$~is a nontrivial Banach space, and the operator $A:\XX\to H$ taking $(v,\eta)$ to~$v(0)$ is surjective. Thus, by Theorem~\ref{T:Jmin-A}, Problem~\ref{probl}
has a unique minimizer $(\bar v^{\hat u,\tau},\,\bar\eta^{\hat u,\tau})$, which linearly depends on~$w_0$. Inequality~\eqref{22} follows immediately from~\eqref{14}, because the norm of~$\bar\eta^{\hat u,\tau}(w_0)$ in the space $L^2(I_{\tau},\,E_M)$ is necessarily smaller than the norm of the control function~$\eta$ constructed in Lemma~\ref{L:exist}.
\end{proof}

{\em Proof of Theorem \ref{Th:main.lin}}. 
Let us fix a sufficiently large $N\ge1$ and denote by~$M_1$ the integer $M$ constructed in Lemma~\ref{L:exist}. 
The main idea of the proof was outlined in the introduction: we use the operator~$\bar\eta^{\hat u,\tau}$
constructed in Lemma~\ref{L:min-lin.cont} to define an exponentially stabilizing control~$\eta^{\hat u,\lambda}$ consecutively
on the intervals~$I_n=(n,n+1)$, $n\ge0$. Namely, let us fix an initial function $v_0\in H$ and set\,\footnote{Recall that the operator $\bar\eta^{\hat u,\tau}$ depends on~$N$.}
$$
\eta^{\hat u,\lambda}(t)=\bar\eta^{\hat u,0}(v_0)(t)\quad\mbox{for} \quad t\in I_0.
$$
Assuming that $\eta^{\hat u,\lambda}$ is constructed on the interval~$(0,n)$
and denoting by~$v(t)$ the corresponding solution on~$[0,n]$, we define
$$
\eta^{\hat u,\lambda}(t)=\bar\eta^{\hat u,n}(v(n))(t)\quad\text{for}\quad t\in I_n.
$$
By construction, $\eta^{\hat u,\lambda}$ is an $E_{M_1}$-valued function square integrable on every bounded interval. Moreover,
the linearity of~$\bar\eta^{\hat u,\tau}$ implies that~$\eta^{\hat u,\lambda}$ linearly depends on~$v_0$.
We claim that, if $N\in\mathbb N$ is sufficiently large, then the solution~$v$ of system~\eqref{eq.v}, \eqref{v0}
with~$\eta=\eta^{\hat u,\lambda}$ satisfies inequalities~\eqref{ed1} and~\eqref{ed2}.

Indeed, it follows from \eqref{l2} that
$$
|v(1)|_V^2=|S_{\hat u,0}(v_0,\,\bar\eta^{\hat u,0}(v_0))(1)|_V^2
\leq \overline C_{[|\hat u|_{\WW}]}\left(|v_0|_H^2
+|\chi|_{L^\infty(\Omega)}^2 |\bar\eta^{\hat u,0}(v_0)|_{L^2(I_0,\,E_M)}^2\right),
$$
where we set $E=E_{M_1}$ to simplify the notation.
Since $\Pi_Nv(1)=0$, we obtain
$$
\alpha_N|v(1)|_H^2\leq |v(1)|_V^2\leq \overline C_{[|\hat
u|_{\mathcal W}]}\left(|v_0|^2+|\chi|_{L^\infty(\Omega)}^2|\bar\eta^{\hat
u,0}(v_0)|_{L^2(I_0,\,E_M)}^2\right).
$$
Using the continuity of $\bar\eta^{\hat u,0}$ (see Lemma~\ref{L:min-lin.cont})
and setting $C_\chi^\prime:=C_\chi|\chi|_{L^\infty(\Omega)}^2$, we derive
$$
|v(1)|_H^2\leq\alpha_N^{-1}|v(1)|_V^2\leq\alpha_N^{-1}(\overline
C_{[|\hat u|_{\mathcal W}]} +C_\chi^\prime)|v_0|^2.
$$
Taking $N$ so large that $\alpha_N\geq e^{\lambda}(\overline
C_{[|\hat u|_{\mathcal W}]}+C_\chi^\prime)$, we obtain
$$
|v(1)|_H^2\leq e^{-\lambda}|v_0|_H^2.
$$
We may repeat the above argument on every the interval~$I_n$ and
conclude that
$$
|v(n+1)|_H^2\leq e^{-\lambda}|v(n)|_H^2.
$$
By induction, we see that the solution~$v$ of problem~\eqref{eq.v}, \eqref{v0} with~$\eta=\eta^{\hat u,\lambda}$
satisfies the inequality
\begin{equation}\label{exp.int}
|v(n)|_H^2\leq e^{-\lambda n}|v_0|_H^2.
\end{equation}
On the other hand, in view of~\eqref{l2}, we have
$$
|v|_{C(\bar I_n,\,H)}^2
\leq \overline C_{[|\hat u|_{\mathcal W}]}\bigl(|v(n)|_H^2
+|\chi|_{L^\infty}^2|\bar\eta^{\hat u,n}(v(n))|_{L^2(I_n,\,E_M)}^2\bigr)
\leq (\overline C_{[|\hat u|_{\mathcal W}]} +C_\chi^\prime)|v(n)|_H^2.
$$
Combining this with~\eqref{l1}, we see that~$v$ satisfies inequality~\eqref{ed1}.

\smallskip
We now prove~\eqref{3.4}. 
It follows from~\eqref{22} and~\eqref{exp.int} that, for any $\tilde\lambda<\lambda$, we have
\begin{align}
|e^{(\tilde\lambda/2)t}\eta^{\hat u,\lambda}|_{L^2(\R_+,\,E_M)}^2
&=\sum_{n\in\mathbb N}|e^{(\tilde\lambda/2)t}
\bar\eta^{\hat u,n}(v(n))|_{L^2(I_n,\,E_M)}^2
\leq C_\chi^\prime\sum_{n\in\mathbb N}e^{\tilde\lambda(n+1)}|v(n)|_H^2\notag\\
&\leq C_\chi^\prime e^{\tilde\lambda}\sum_{n\in\mathbb N} e^{(\tilde\lambda-\lambda)n}|v_0|_H^2
\leq \kappa_2|v_0|_H^2.
\label{bdd.st-control}
\end{align}

It remains to prove inequality~\eqref{ed2}.
In view of~\eqref{l2}, we have
\begin{align*}
|\sqrt{t-n} v|_{C(\bar I_n,\,V)}^2 
&\leq \overline C_{[|\hat u|_{\mathcal W}]}\left(|v(n)|_H^2
+3|\chi|_{L^\infty(\Omega)}^2|\bar\eta^{\hat u,n}(v(n))|_{L^2(I_n,\,E_M)}^2\right) \\
&\leq \bigl(\overline C_{[|\hat u|_{\mathcal W}]} +3C_\chi^\prime\bigr)|v(n)|_H^2.
\end{align*}
Combining this with~\eqref{exp.int}, we see that
$$
|v(n+1)|_V^2\leq C_1 |v(n)|_H^2\leq C_1e^{-\lambda n}|v_0|_H^2.
$$
In view of~\eqref{l3}, for $n\ge1$ we derive
\begin{align*}
|v|_{C(\bar I_n,\,V)}^2\leq
\overline C_{[|\hat u|_{\mathcal W}]}|v(n)|_V^2
+3|\chi|_{L^\infty(\Omega)}^2|\bar\eta^{\hat u,n}(v(n))|_{L^2(I_n,\,E_M)}^2
\leq C_2 e^{-\lambda(n-1)}|v_0|_H^2,
\end{align*}
whence it follows that
$$
|v(t)|_V^2\leq C_3e^{-\lambda t}|v_0|_H^2\quad\mbox{for}\quad t\ge1.
$$
Using again inequality~\eqref{l3}, we conclude that~\eqref{ed2} holds. The proof of the theorem is complete.
\endproof

\subsection{Feedback control}
In this section, we show that the exponentially stabilizing control constructed in Theorem~\ref{Th:main.lin} can be chosen in a feedback form.
Namely, let us fix a nonzero function $\chi\in C_0^\infty(\Omega)$ and denote by~$\EE_M$ the vector
space spanned by the functions $\chi\phi_j$, $j=1,\dots,M$. Note that, due to elliptic regularity, the eigenfunctions~$\phi_j$ are infinitely smooth in~$\Omega$, whence it follows, in particular, that~$\EE_M$ is contained in $C_0^\infty(\omega,\R^3)$ for any $\omega\supset\supp\chi$.  We will prove the following theorem.

\begin{theorem} \label{feedback}
Given $\hat u\in\WW$ and $\lambda>0$, let $M=\overline C_{[|\hat u|_\WW,\lambda]}\in\mathbb N$ be the integer constructed in Theorem~\ref{Th:main.lin}. Then there are a family of continuous operators $K_{\hat u}^\lambda (t):H\to\EE_M$ and a constant $\kappa=\overline C_{[|\hat u|_\WW,\lambda]}$ such that
the following properties hold.
\begin{itemize}
 \item[\bf(i)]
The function $t\mapsto K_{\hat u}^\lambda (t)$ is continuous in the weak operator topology, and its operator norm is bounded
by~$\kappa$.
 \item[\bf(ii)]
For any $s\ge0$ and $v_0\in H$, the solution of the problem
\begin{align}
v_t + Lv+{\mathbb B}(\hat u)v &=\Pi K_{\hat u}^\lambda (t)v,\label{eq.vk}\\
v(s)&=v_0\label{vs}
\end{align}
exists on the half-line~$\R_s$ and satisfies the inequality
\begin{equation} \label{ed1-feedback}
e^{\lambda(t-s)}|v(t)|_H^2+\int_s^t e^{\lambda (\tau-s)}\bigl(|v(\tau)|_V^2+|v_t(\tau)|_{V'}^2\bigr)d\tau
\leq \kappa|v_0|_H^2, \quad t\ge s,
\end{equation}
Moreover, if $v_0\in V$, then
\begin{equation} \label{ed2-feedback}
e^{\lambda(t-s)}|v(t)|_V^2+\int_s^t e^{\lambda (\tau-s)}\bigl(|v(\tau)|_{{\mathrm D}(L)}^2+|v_t(\tau)|_{H}^2\bigr)d\tau
\leq \kappa |v_0|_V^2, \quad t\ge s.
\end{equation}
\end{itemize}
\end{theorem}

To prove this theorem, we will need two auxiliary lemmas. Let us consider the following problem.

\begin{problem}\label{probl.Ms}
Given $s\ge0$, $\lambda>0$, $M\in\mathbb N$ and $w_0\in H$, find the minimum of the functional
$$
M_s^\lambda(v,\,\eta):=\int_{\R_s}e^{\lambda t}(|v|_V^2+|\eta|_{L^2}^2)\,{\rm d}t
$$
on the set of functions $(v,\eta)\in W_\lambda(\R_s,\,V,\,V^\prime)\times L_\lambda^2(\R_s,\,L^2(\Omega,\R^3))$
that satisfy equation~\eqref{eq.v} and the initial condition
\begin{equation} \label{ICs}
v(s)=w_0.
\end{equation}
\end{problem}

The following lemma establishes the existence of an optimal solution and gives a formula for the optimal cost.

\begin{lemma} \label{Ms}
For any $\hat u\in\mathcal W$ and $\lambda>0$ there is an integer
$M=\overline C_{[|\hat u|_{\mathcal W},\lambda]}\ge1$ such that
Problem~\ref{probl.Ms} has a unique minimizer~$(v_s^*,\eta_s^*)$.
Moreover, there is a continuous operator $Q_{\hat
u}^{s,\lambda}:H\to H$ such that
\begin{align}
M_s^\lambda(v_s^*,\,\eta_s^*)&=\bigl(Q_{\hat u}^{s,\lambda}w_0,w_0\bigr),\label{optimalcost}\\
|Q_{\hat u}^{s,\lambda}|_{\LL(H)}&\leq Ce^{\lambda s},\label{normofQ}
\end{align}
where $\LL(H)$ stands for the space of continuous linear operators in~$H$ with the natural norm and $C=\overline C_{[|\hat u|_{\mathcal W},\lambda]}>0$ is a
constant. Finally, $Q_{\hat u}^{s,\lambda}$ continuously depends
on~$s$ in the weak operator topology.
\end{lemma}

\begin{proof}
Let~$\mathcal X$ be the space of functions $(v,\eta)\in
W_\lambda(\R_s,\,V,\,V^\prime)\times
L_\lambda^2(\R_s,\,L^2)$ that
satisfy~\eqref{eq.v} and endow it with the norm
$M_s^\lambda(v,\,\eta)^{1/2}$. It is straightforward to see
that~$\XX$ is a Hilbert space. Using
Theorem~\ref{Th:main.lin} with a constant $\hat\lambda>\lambda$
and the initial point moved to~$s$, one can construct an integer
$M=\overline C_{[|\hat u|_{\WW},\hat\lambda]}\geq1$ such that, for
any $w_0\in H$ and an appropriate control~$\eta\in L_{\hat\lambda}^2(\R_s,E_M)$, we have
$$
|v(t)|_H^2\leq \kappa_1 e^{-\hat\lambda(t-s)}|w_0|_H^2,\quad 
|\eta(t)|_{E_M}^2\leq \kappa_2 e^{-\hat\lambda(t-s)}|w_0|_H^2,\quad
$$
where $v$ stands for the solution of~\eqref{eq.v}, \eqref{ICs}. Furthermore, by Lemma~\ref{bd.hatu-r0}, we have
$$
\int_{I_\tau}e^{\lambda t}|v|_V^2{\rm d}t
\le e^{\lambda(\tau+1)}\int_{I_\tau}|v|_V^2{\rm d}t
\le {\overline C}_{[|\hat u|_\WW,\lambda]}\,e^{\lambda \tau}|v(\tau)|_H^2
\quad\mbox{for any $\tau\ge0$}. 
$$
Combining the above three inequalities, we conclude that 
\begin{align}
M_s^\lambda(v,\,\eta)
&=\int_{\R_s}e^{\lambda t}(|v|_V^2+|\eta|_{L^2}^2)\,{\rm d}t
\leq \overline C_{[\hat\lambda,(\hat\lambda-\lambda)^{-1},|\hat u|_{\mathcal W}]} 
\,e^{\lambda s}|w_0|_H^2.\label{est.Ms}
\end{align}
It follows that~$\XX$ is nonempty, and the mapping $A:\XX\to H$ taking~$(v,\eta)$ to~$v(0)$ is surjective. Thus, by Theorem~\ref{T:Jmin-A},
Problem~\ref{probl.Ms} has a unique minimizer~$(v_s^*,\eta_s^*)=(v_s^*(w_0),\eta_s^*(w_0))$, which linearly depends on~$w_0$.

\smallskip
We now prove~\eqref{optimalcost} and~\eqref{normofQ}. It follows from~\eqref{est.Ms} that the mapping
$$
(a,\,b)\mapsto \int_{\R_s}e^{\lambda
t}((v_s^\ast(a),\,v_s^\ast(b))_V+(\eta_s^\ast(a),\,\eta_s^\ast(b))_{L^2})\,{\rm d}t
$$
is a continuous bilinear form on~$H$ which is bounded by $C_2e^{\lambda s}$ on the unit ball. Therefore, the optimal cost
can be written as~\eqref{optimalcost}, where~$Q_{\hat u}^{s,\lambda}$ is a bounded self-adjoint operator in~$H$ whose norm satisfies~\eqref{normofQ}.

\smallskip
It remains to establish the continuity of~$Q_{\hat u}^{s,\lambda}$ in the weak operator topology. To this end, it suffices to prove that
\begin{equation} \label{3.26}
(Q_{\hat u}^{s,\lambda}w,w)\to(Q_{\hat u}^{s_0,\lambda}w,w)
\quad\mbox{as $s\to s_0$ for any $w\in H$}.
\end{equation}
We will prove this convergence for $s_0>0$; in the case $s_0=0$, the proof is simpler.

Let us fix $\tau\le s_0$ and denote by $(v_\tau^*,\eta_\tau^*)$ the optimal solution of Problem~\ref{probl.Ms} with $s=\tau$ and $w_0=w$. Note that, for any bounded interval $I\subset\R_\tau$, the norms of the functions~$v_\tau^*$ and~$\eta_\tau^*$ in the spaces $W(I,V,V')$ and~$L^2(I,L^2)$, respectively, are bounded uniformly in~$\tau$. It is straightforward to see that the restriction of~$(v_\tau^*,\eta_\tau^*)$ to the half-line~$\R_s$ with $s>\tau$ is the optimal solution of Problem~\ref{probl.Ms} with $w_0=v_\tau^*(s)$; cf.\ Lemma~\ref{Ns} below. Therefore, abbreviating $Q^s=Q_{\hat u}^{s,\lambda}$, we can write
$$
\bigl(Q^sv_\tau^*(s),v_\tau^*(s)\bigr)
=\int_s^\infty e^{\lambda t}\bigl(|v_\tau^*(t)|_V^2+|\eta_\tau^*(t)|_{L^2}^2\bigr)\,{\rm d}t.
$$
Setting $\Delta_\tau^s(w)=|(Q^sv_\tau^*(s),v_\tau^*(s))-(Q^sw,w)|$,  for $s\ge\tau$ we have
\begin{equation*} 
\bigl|(Q^sw,w)-(Q^{s_0}w,w)\bigr|\le\Delta_\tau^s(w)+\Delta_\tau^{s_0}(w)
+\biggl|\int_{s_0}^s e^{\lambda t}\bigl(|v_\tau^*(t)|_V^2+|\eta_\tau^*(t)|_{L^2}^2\bigr)\,{\rm d}t\biggr|.
\end{equation*}
The third term on the right-hand side of this inequality goes to zero as $s\to s_0$. Therefore convergence~\eqref{3.26} will be established if we prove that $\Delta_\tau^s(w)\to0$ as $\tau,s\to s_0$. To this end, note that the above-mentioned boundedness of~$v_\tau^*$ implies that 
$$
|v_r^*(s)-w|_{L^2}\to0 \quad\mbox{as $\tau,s\to s_0$}.
$$
Combining this with the fact that the norm of $Q^s$ is bounded on finite intervals (see~\eqref{normofQ}), we arrive at the required assertion. 
\end{proof}

We now consider another minimization problem closely related to Problem~\ref{probl.Ms} with $s=0$.

\begin{problem}\label{probl.Ns}
Given $\lambda>0$ and $v_0\in H$, find the minimum of the functional
$$
N_s^\lambda(v,\,\eta):=\int_{(0,\,s)}e^{\lambda t}(|v|_V^2
+|\eta|_{L^2}^2)\,{\rm d}t+(Q_{\hat u}^{s,\lambda} v(s),\,v(s))
$$
on the set of functions $(v,\,\eta)\in
W((0,\,s),\,V,\,V^\prime)\times L^2((0,\,s),\,L^2(\Omega,\R^3))$ that
satisfy~\eqref{eq.v}, \eqref{v0}, where~$M$ is the integer
constructed in Lemma~\ref{Ms}.
\end{problem}

Theorem~\ref{T:Jmin-A} implies that Problem~\ref{probl.Ns} has a unique minimizer $(v_s^\bullet,\,\eta_s^\bullet)$,
which is a linear function of~$v_0\in H$. The following lemma is the dynamic programming principle for Problem~\ref{probl.Ms} with~$s=0$.

\begin{lemma} \label{Ns}
Under the hypotheses of Lemma~\ref{Ms}, the restriction of~$(v_0^*,\eta_0^*)$ to the interval $(0,s)$
coincides with $(v_s^\bullet,\,\eta_s^\bullet)$ and the restriction of~$(v_0^*,\eta_0^*)$ to the half-line~$\R_s$
coincides with $(v_s^*,\eta_s^*)(v_0^*(s))$.
\end{lemma}

\begin{proof}
We will confine ourselves to the proof of the first assertion,
because the second one is obvious. Let us define the function
$$
(z_0^\ast,\,\theta_0^\ast)(t):=\begin{cases}
                                 (v_s^\bullet,\,\eta_s^\bullet)(v_0)(t)
                                 &\quad\text{for}\quad t\in(0,\,s),\\
                 (v_s^\ast,\,\eta_s^\ast)(v_s^\bullet(s))(t)
                 &\quad\text{for}\quad t\in\R_s.
                               \end{cases}
$$
Then we have
$$
M_0^\lambda(z_0^\ast,\,\theta_0^\ast)=N_s^\lambda(v_s^\bullet,\,\eta_s^\bullet).
$$
On the other hand, the definition of $(v_s^\bullet,\,\eta_s^\bullet)$ implies that
$$
N_s^\lambda(v_s^\bullet,\,\eta_s^\bullet)
\leq N_s^\lambda\bigl((v_0^\ast,\,\eta_0^\ast)|_{(0,s)}\bigr)
\leq M_0^\lambda(v_0^\ast,\,\eta_0^\ast),
$$
whence it follows that
$$
M_0^\lambda(z_0^\ast,\,\theta_0^\ast) \leq M_0^\lambda(v_0^\ast,\,\eta_0^\ast).
$$
The uniqueness of minimizer for Problem~\ref{probl.Ms} with $s=0$ implies that
$(z_0^\ast,\,\theta_0^\ast)=(v_0^\ast,\,\eta_0^\ast)$, and the required assertion follows.
\end{proof}

{\em Proof of Theorem~\ref{feedback}}. 
{\it Step~1}. It is straightforward to see that Problem~\ref{probl.Ns}
satisfies the hypotheses of the Karush--Kuhn--Tucker
Theorem~\ref{T:Jmin-F}, in which
$$\XX=W((0,\,s),V,V')\times L^2((0,\,s),\,L^2(\Omega,\R^3)), \quad \YY=H\times L^2((0,\,s),\,V),
$$
$J=N_s^\lambda$, and $F:\XX\to\YY$ is the affine operator taking~$(v,\eta)$ to $\bigl(v(0)-w_0, v_t+Lv+{\mathbb B}(\hat u)v-\Pi(\chi P_M\eta)\bigr)$.
Hence, there is a Lagrange multiplier $(\mu_s,\, q_s)\in H\times L^2((0,\,s),\,V)$ such that
$$
(N_s^\lambda)^\prime(v_s^\bullet,\,\eta_s^\bullet)-(\mu_s,\,q_s)
\circ F^\prime(v_s^\bullet,\,\eta_s^\bullet)=0.
$$
It follows that, for all $z\in W((0,\,s),\,V,\,V^\prime)$ and $\xi\in
L^2((0,\,s),\,L^2)$, we have
\begin{align}
&2\int_0^se^{\lambda t}(v_s^\bullet,\,z)_V\,{\rm d}t+2(Q_{\hat u}^{s,\lambda}v_s^\bullet(s),\,z(s))_H+(z(0),\,\mu_s)_H\notag\\
&\qquad +\int_0^s\langle z_t+Lz
+\mathbb B(\hat u)z,\, q_s\rangle_{V^\prime,V}\,{\rm d}t=0,
\label{s.z}\\
&2\int_0^se^{\lambda t}(\eta_s^\bullet,\,\xi)_{L^2}\,{\rm d}t+\int_0^s\langle -\Pi(\chi P_M\xi),\,
q_s\rangle_{V^\prime,V}\,{\rm d}t=0.\label{s.xi}
\end{align}
In particular, we conclude from~\eqref{s.z}  that
\begin{equation}\label{eq.qs}
(q_s)_t-Lq_s-\mathbb B^\ast(\hat u)q_s=2e^{\lambda t}Lv_s^\bullet(t).
\end{equation}
Since $q_s,v_s^\bullet\in L^2((0,s),V)$, we see that $\p_tq_s\in
L^2((0,s),V')$, and therefore $q_s\in W((0,\,s),V,\,V^\prime)$,
whence it follows that $q_s\in C([0,\,s],\,H)$. Using
again~\eqref{s.z}, we derive
\begin{equation} \label{qs1}
q_s(s)=-2Q_{\hat u}^{s,\lambda}v_s^\bullet(s).
\end{equation}
On the other hand, relation~\eqref{s.xi} implies that
\begin{equation}\label{etas.qs}
\eta_s^\bullet=\frac{1}{2}e^{-\lambda t}P_M(\chi q_s).
\end{equation}
In particular, $\eta_s^\bullet(t)$ is a continuous function of $t$ with range in~$E_M$.
Combining~\eqref{etas.qs} and~\eqref{qs1}, we derive
$$
\eta_s^\bullet(s)=-e^{-\lambda s}P_M\chi Q_{\hat u}^{s,\lambda}v_s^\bullet(s).
$$
Recalling Lemma~\ref{Ns} and using the fact that~$s$ is arbitrary, we conclude that~$\eta_0^*$ is a continuous function of
time with range in~$E_M$ and that
$$
\eta_0^\ast(t)=-e^{-\lambda t}P_M\bigl(\chi Q_{\hat u}^{t,\lambda}v_0^\ast(t)\bigr)\quad\mbox{for all}\quad t\ge0.
$$
Thus, the optimal trajectory $v_0^\ast$ for Problem~\ref{probl.Ms} with $s=0$ satisfies~\eqref{eq.vk}, 
where 
$$
K_{\hat u}^\lambda(t):=-e^{-\lambda t}\chi P_M\chi Q_{\hat u}^{t,\lambda}.
$$
It is clear that $K_{\hat u}^\lambda(t)$ is a linear continuous
operator from~$H$ to~$\EE_M$. Moreover, it continuously depends
on~$t$ in the weak operator topology, because so does the
family~$Q_{\hat u}^{t,\lambda}$. Finally, it follows
from~\eqref{normofQ} that the norm of~$K_{\hat u}^\lambda(t)$ is
bounded by a constant depending only on~$\lambda$ and~$|\hat
u|_\WW$. We have thus constructed a feedback control~$K_{\hat
u}^\lambda(t)$ possessing the properties mentioned in~(i). Moreover, repeating the above arguments for Problem~\ref{probl.Ms} with an arbitrary~$s>0$, we conclude that its optimal solution $(v_s^*,\eta_s^*)$ satisfies the relation $\eta_s^*(t)=-e^{-\lambda t}P_M(\chi Q_{\hat u}^{t,\lambda}v_s^*(t))$ for $t\ge s$. Hence, if $v(t)$ is the solution of problem~\eqref{eq.vk}, \eqref{vs}, then 
$$
\bigl(v(t),-e^{-\lambda t}P_M(\chi Q_{\hat u}^{t,\lambda}v(t))\bigr)
=\bigl(v_s^*(t),\eta_s^*(t)\bigr)\quad\mbox{for $t\ge s$}. 
$$
Combining this with~\eqref{normofQ}, we conclude that
\begin{equation} \label{3.22}
\bigl(Q_{\hat u}^{s,\lambda}v_0,v_0\bigr)
=\int_{\R_s}\bigl(e^{\lambda t}|v(t)|_V^2
+e^{-\lambda t}\bigl|P_M(\chi Q_{\hat u}^{t,\lambda}v(t))\bigr|_{L^2}^2\bigr)\,{\rm d}t
\le C\,e^{\lambda s}|v_0|_H^2.
\end{equation}

{\it Step~2}. We now prove inequalities~\eqref{ed1-feedback} and~\eqref{ed2-feedback} for solutions of problem~\eqref{eq.vk}, \eqref{vs}.
Let us fix~$v_0\in H$ and denote by~$v$ the solution of~\eqref{eq.vk}, \eqref{vs}. It is straightforward to see that the function $z(t)=e^{(\lambda/2) t}v(t)$ satisfies the equation
\begin{equation} \label{exp.v0*}
z_t+Lz+\mathbb B(\hat u)z=\frac{\lambda}{2}z+K_{\hat u}^\lambda(t)z.
\end{equation}
Taking the scalar product of~\eqref{exp.v0*} with $2z$ and using the uniform boundedness of the family~$K_{\hat u}^\lambda(t)$, we derive
\begin{equation*}
\frac{\rm d}{{\rm d} t}|z(t)|_H^2+2|z(t)|_V^2 =\lambda |z(t)|_H^2
+2(K_{\hat u}^\lambda(t)z,z) -2(\mathbb B(\hat u)z(t),\,z(t))_H
\leq {\overline C}_{[|\hat u|_\WW,\lambda]}|z(t)|_H^2
\end{equation*}
Integrating this inequality over the interval $(s,\,t)$ with $t>s$, recalling the definition of~$z$, and using Lemmas~\ref{Ms}, \ref{Ns} and inequality~\eqref{3.22}, we obtain
\begin{align}
e^{\lambda t}|v(t)|_H^2+2\int_{(s,\,t)}e^{\lambda\tau}|v(\tau)|_V^2 \,{\rm d}\tau
&\leq e^{\lambda s}|v_0|_H^2
+C_1\int_{(s,\,t)}e^{\lambda \tau}|v(\tau)|_H^2\,{\rm d}\tau\notag\\
&\leq e^{\lambda s}|v_0|_H^2
+C_1\int_{\R_s}e^{\lambda \tau}|v(\tau)|_H^2\,{\rm d}\tau\notag\\
&\leq C_2 e^{\lambda s}|v_0|_H^2. \label{ed3}
\end{align}
Furthermore, it follows from~\eqref{eq.vk} that
$$
e^{\lambda t}|v_t(t)|_{V'}^2\le C_3 e^{\lambda t}|v(t)|_V^2.
$$
Combining this with~\eqref{ed3}, we arrive at the inequality
$$
e^{\lambda
t}|v(t)|_H^2+2\int_{(s,\,t)}e^{\lambda\tau}\bigl(|v(\tau)|_V^2+|v_t(\tau)|_{V'}^2\bigr)\,{\rm
d}\tau \leq C_4e^{\lambda s}|v_0|_H^2,
$$
which is equivalent to~\eqref{ed1-feedback}.

\smallskip
We now assume that $v_0\in V$. Taking the scalar product of~\eqref{exp.v0*} with $2Lz$, and using the Schwarz inequality and
the uniform boundedness of the family~$K_{\hat u}^\lambda(t)$, we derive
\begin{align*}
\frac{\rm d}{{\rm d} t}|z(t)|_V^2+2|z(t)|_{{\rm D}(L)}^2
&=\lambda(z,Lz)_H+2(K_{\hat u}^\lambda(t)z,Lz)-2(\mathbb B(\hat u)z(t),\,Lz(t))_H \\
&\le |Lz|_H^2+{\overline C}_{[|\hat u|_\WW,\lambda]} |z(t)|_V^2.
\end{align*}
Integrating this inequality over the interval $(s,\,t)$ and using~\eqref{ed1-feedback}, we obtain
\begin{align} 
e^{\lambda t}|v(t)|_V^2+\int_{(s,\,t)}e^{\lambda\tau}|v(\tau)|_{{\rm D}(L)}^2 \,{\rm d}\tau
&\leq e^{\lambda s}|v_0|_V^2+C_5\int_{(s,t)}
e^{\lambda \tau}|v(\tau)|_V^2\,{\rm d}\tau\notag\\
&\le C_6e^{\lambda s}|v_0|_V^2.\label{48}
\end{align}
Furthermore, relation~\eqref{eq.vk} implies that
$$
e^{\lambda t}|v_t(t)|_{H}^2\le C_6 e^{\lambda t}|v(t)|_{{\rm
D}(L)}^2.
$$
Combining this with~\eqref{48}, we arrive at~\eqref{ed2-feedback}.
\endproof

\bigskip
We conclude this section with a few remarks. 
\begin{remarks}
\rm 
{\bf(a)}
Once a feedback control is constructed, it is easy to find a time-dependent Lyapunov function for the problem in question. Indeed, let $U(s,t)$ be the operator  taking~$v_0\in H$ to~$v(t)$, where~$v$ stands for the solution of~\eqref{eq.vk}, \eqref{vs}. It is straightforward to check that the functional 
$$
\varPhi(t,w)=\int_t^\infty |U(t,\tau)w|_{L^2(\Omega,\R^3)}^2{\rm d}\tau
$$
decays along the trajectories of~\eqref{eq.vk}. It is difficult, however, to write down this functional in a more explicit form. 

\smallskip
{\bf(b)}
The operator $Q_{\hat u}^{t,\lambda}$ defining the optimal cost satisfies the following Ricatti equation:
\begin{equation} \label{riccati}
\dot Q-(Q\,{\mathbb L}(\hat u)+{\mathbb L}^*(\hat u)Q)
-e^{-\lambda s}Q(\Pi\chi P_M\chi\Pi)Q=-e^{\lambda s}L, 
\end{equation}
where ${\mathbb L}(\hat u)=L+{\mathbb B}(\hat u)$. 
Since this equation does not play any role in this paper, we confine ourselves to its formal derivation. Let~$v$ be the solution of~\eqref{eq.vk}, \eqref{vs} with $v_0\in H$ and let $\eta(t)=-e^{-\lambda t}P_M\chi Q^tv(t)$, where we set $Q^t=Q_{\hat u}^{t,\lambda}$. By the dynamic programming principle (cf.\ Lemma~\ref{Ns}), the restriction of~$(v,\eta)$ to the half-line~$\R_\tau$ is the optimal solution of Problem~\ref{probl.Ms} with $s=\tau$ and $w_0=v(\tau)$. Therefore, we have
(cf.\ the equality in~\eqref{3.22}) 
$$
\bigl(Q^\tau v(\tau),v(\tau)\bigr)
=\int_{\R_\tau}\bigl(e^{\lambda t}|v(t)|_V^2
+e^{-\lambda t}\bigl|P_M(\chi Q^tv(t))\bigr|_{L^2}^2\bigr)\,{\rm d}t.
$$
Differentiating this relation with respect to~$\tau$ and carrying out some simple transformations, we obtain
$$
\bigl(\bigl(\dot Q^\tau-(Q^\tau{\mathbb L}(\hat u)+{\mathbb L}^*(\hat u)Q^\tau)
-e^{-\lambda \tau}Q^\tau(\Pi\chi P_M\chi\Pi)Q^\tau+e^{\lambda\tau}L\bigr)v(\tau),v(\tau)\bigr)=0. 
$$
Setting $\tau=s$ and recalling that $v(s)=v_0$ is arbitrary, we conclude that~$Q^s$ satisfies~\eqref{riccati}. 

\smallskip
{\bf(c)}
In the case of a stationary reference solution~$\hat u$, it is possible to give a rather  sharp description of the dimension~$M$ for the feedback controller whose range depends on~$\hat u$; e.g., see~\cite{BT-2004,BT-2009,RT-2009}. In our situation, the range of the controller depends only on the norm of~$\hat u$, and its space dimension is determined by the integer~$M_1$ in the truncated observability inequality~\eqref{obs-in.M} with a sufficiently large integer~$N$. However, the feedback operator depends on time, and its image may be infinite-dimensional in time. It would be interesting to find out if it is possible to reduce the space dimension of the controller in our situation using further information about~$\hat u$. 
\end{remarks}

\section{Stabilization of the nonlinear problem}\nopagebreak
\label{S:nonlinear}
\subsection{Main result}
Let us consider the nonlinear problem
\begin{align}
v_t + Lv+Bv+\mathbb B(\hat u)v
&=K_{\hat u}^\lambda(t)v,\quad t\in \R_+;\label{eq.tildev}\\
v(0)&=v_0,\label{tildev0}
\end{align}
where the operator $K_{\hat u}^\lambda(t)$ is constructed in
Theorem~\ref{feedback}. Given a constant~$\lambda>0$, we denote
by~$\mathcal Z^\lambda$ the space of functions $z\in
C(\R_+,\,V)\cap L^2_{\rm loc}(\R_+,\,U)$
 such that
$$
|z|_{\mathcal Z^\lambda}
:=\sup_{t\geq 0}\left(e^{\lambda t}|z(t)|_V^2+\int_{(t,\,t+1)}e^{\lambda \tau}|z(\tau)|_{{\rm D}(L)}^2\,{\rm d}\tau\right)^{1/2}<\infty.
$$
The following theorem is the main result of this paper.

\begin{theorem}\label{Th:main}
Given $\hat u\in \mathcal W$ and~$\lambda>0$, let $M=\overline C_{[|\hat u|_\WW,\lambda]}$ be
the integer constructed in Theorem~\ref{feedback}. Then there are
positive constants~$\vartheta$ and~$\epsilon$ depending only on
$|\hat u|_{\mathcal W}$ and~$\lambda$ such that for
$|v_0|_V\le\epsilon$ the solution~$v$ of system \eqref{eq.tildev},
\eqref{tildev0} is well defined for all $t\ge0$ and satisfies the
inequality
\begin{equation} \label{ed}
|v(t)|_V^2\leq\vartheta e^{-\lambda t}|v_0|_V^2\quad\mbox{for $t\ge0$}.
\end{equation}
\end{theorem}

\begin{proof}
We will use the contraction mapping principle. We fix a constant~$\vartheta>0$ and a function
$v_0\in V$ and introduce the following subset of~$\mathcal Z^\lambda$:
$$
\mathcal Z^\lambda_\vartheta:=\{z\in\mathcal Z^\lambda\mid\,z(0)=v_0,
|z|_{\mathcal Z^\lambda}^2\leq \vartheta|v_0|_V^2\}.
$$
We define a mapping $\Xi:\mathcal Z^\lambda_\vartheta\to
C(\R_+,V)\cap L^2_{\rm loc}(\R_+,\,U)$ that takes a function $a\in
\mathcal Z^\lambda$ to the solution of the problem
\begin{align}
b_t + Lb+\mathbb B(\hat u)b&=K_{\hat u}^\lambda b-Ba,\quad t\in \R_+,\label{eqb}\\
b(0)&=v_0.\label{icb}
\end{align}
Suppose we have shown the following proposition.

\begin{proposition}\label{L:contract}
Under the hypotheses of Theorem~\ref{Th:main}, there exists $\vartheta>0$ such that the following property holds: for any $\gamma\in(0,\,1)$  {\color{blue}{one can find a}} constant $\epsilon=\epsilon_\gamma>0$ {\color{blue}{such that for any $v_0\in V$ with $|v_0|_V\le\epsilon$}} the mapping~$\Xi$
takes the set~$Z^\lambda_\vartheta$ into itself and satisfies the inequality
\begin{equation} \label{contraction}
|\Xi(a_1)-\Xi(a_2)|_{\mathcal Z^\lambda}\leq\gamma|a_1-a_2|_{\mathcal Z^\lambda}\quad
\mbox{for all}\quad a_1,\,a_2\in\mathcal Z^\lambda_\vartheta.
\end{equation}
\end{proposition}

Thus, if $|v_0|_V$ is sufficiently small, then the contraction
mapping principle implies that there is a unique fixed point $v\in
\mathcal Z^\lambda_\vartheta$ for~$\Xi$. It follows from the
definition of~$\Xi$ and~$Z^\lambda_\vartheta$ that~$v$ is a
solution of problem~\eqref{eq.tildev}, \eqref{tildev0} and
satisfies the required inequality~\eqref{ed}. We claim that~$v$ is the unique solution of~\eqref{eq.tildev}, \eqref{tildev0} in the space
$C(\R_+,V)\cap L^2_{\rm loc}(\R_+,\,U)$. Indeed, if~$w$ is another solution, then the difference $z=v-w$ vanishes at $t=0$ and satisfies the equation
$$
z_t+Lz+B(z,v)+B(w,z)+\mathbb B(\hat u)z=K_{\hat u}^\lambda(t)z.
$$
Taking the scalar product of this equation with~$z$ in~$H$, carrying out some standard transformations (e.g., see~\cite{temams}), and using the uniform boundedness of the feedback control~$K_{\hat u}^\lambda(t)$ as an operator in~$H$, we see that $z\equiv0$. Hence, to complete the proof of the theorem, it suffices to establish the above proposition. This is done in the next subsection.
\end{proof}

\begin{remark}
\rm 
The hypotheses of Theorem~\ref{Th:main} can be relaxed. Namely, it suffices to assume that the reference solution~$\hat u$ satisfies the condition 
$$
\sup_{\tau\ge0}\bigl(|\hat u|_{L^\infty(Q_\tau)}+|\hat u_t|_{L^2(Q_\tau)}\bigr)<\infty.
$$
Indeed, as is proved in~\cite{FGIP-2004}, the observability inequality~\eqref{obs-in.chi} remains valid in this situation. It follows that the truncated observability inequality~\eqref{obs-in.M}, which is the key point of our approach, is also true. One can check that all the proofs can be carried out under the above weaker hypothesis. However, some calculations become cumbersome, and for the sake of clarity of the paper, we have imposed the more restrictive condition~$\hat u\in\WW$. 
\end{remark}

\subsection{Proof of Proposition~\ref{L:contract}}
{\it Step 1}.
We first derive an estimate for solutions of the equation
\begin{equation} \label{1000}
z_t + Lz+\mathbb B(\hat u)z=K_{\hat u}^\lambda z+f(t),
\end{equation}
where $f\in L_{\rm loc}^2(\R_+,H)$. Namely, we will show that
\begin{equation} \label{est-b}
\sup_{t\ge0}\biggl(e^{\lambda
t}|z(t)|_V^2+\int_{(t,t+1)}e^{\lambda s}|z(s)|_{{\rm D}(L)}^2{\rm
d}s\biggr) \le
C_1\biggl(|z(0)|_V^2+\sup_{t\ge0}\int_{(t,t+1)}e^{2\lambda
s}|f(s)|_H^2{\rm d}s\biggr),
\end{equation}
where $C_1=\overline C_{[|\hat u|_\WW,\lambda]}$ is a constant.
Indeed, recall that $U(s,t)$ denotes  the operator taking~$v_0\in H$
to~$v(t)$, where~$v$ stands for the solution of~\eqref{eq.vk},
\eqref{vs}. By the Duhamel formula, we can write~$z$ as
\begin{equation} \label{duhamel}
z(t)=U(0,\,t)z(0)+\int_{(0,\,t)}U(s,\,t)f(s)\,{\rm d}s.
\end{equation}
Combining this with~\eqref{ed1-feedback}, we derive
\begin{align}
|z(t)|_H^2&=2|U(0,\,t)z(0)|_H^2+2\biggl(\int_{(0,\,t)}|U(s,\,t)f(s)|_H\,{\rm d}s\biggr)^2\notag\\
&=2\kappa e^{-\lambda t}|z(0)|_H^2+2\kappa e^{-\lambda t}
\biggl(\int_{(0,\,t)}e^{(\lambda/2)s}|f(s)|_H\,{\rm
d}s\biggr)^2.\label{101}
\end{align}
Now note that, for any non-negative function $c(t)$ and any $\lambda>0$, we have
\begin{align*}
\sup_{t\ge0}\int_{(0,\,t)}e^{(\lambda/2)s}c(s)\,{\rm d}s &\le
\int_{(0,\,+\infty)}e^{(\lambda/2)s}c(s)\,{\rm d}s=
\sum_{k=1}^\infty \int_{(k-1,\,k)}e^{(\lambda/2)s}c(s)\,{\rm d}s\\
&\le \sum_{k=1}^\infty e^{(\lambda/2)k}\biggl(\int_{(k-1,\,k)}|c(s)|^2\,{\rm d}s\biggr)^{1/2}\\
&\le \sum_{k=1}^\infty e^{-(\lambda/2)(k-2)}\biggl(\int_{(k-1,\,k)}e^{2\lambda s}|c(s)|^2\,{\rm d}s\biggr)^{1/2}\\
&\le C_2 \biggl(\sup_{t\ge0}\int_{(t,\,t+1)}e^{2\lambda
s}|c(s)|^2\,{\rm d}s\biggr)^{1/2}.
\end{align*}
Substituting this inequality with $c(t)=|f(t)|_H$ into~\eqref{101}, we derive
\begin{equation} \label{102}
\sup_{t\ge0}\bigl( e^{\lambda t}|z(t)|_H^2\bigr) \le
2\kappa\biggl(|z(0)|_H^2+C_2^2\sup_{t\ge
0}\int_{(t,\,t+1)}e^{2\lambda s}|f(s)|_H^2\,{\rm d}s\biggr)
\end{equation}
On the other hand, it is easy to see that the analogue of Lemma~\ref{bd.hatu-r0} is true for equation~\eqref{1000}.
In particular, for any $s\ge0$ we have the estimates
\begin{align}
 (t-s)|U(s,t)z_0|_V^2
&\le C_3\biggl(|z_0|_H^2+\int_{(s,t)}|f(\tau)|_H^2{\rm d}\tau\biggr),\label{103}\\
|U(s,t)z_0|_V^2+\int_{(s,s+1)}|U(s,\tau)z_0|_{{\rm D}(L)}^2{\rm d}\tau
&\le C_3\biggl(|z_0|_V^2+\int_{(s,t)}|f(\tau)|_H^2{\rm d}\tau\biggr),\label{104}
\end{align}
where $s\le t\le s+1$, and $C_3>0$ does not depend on~$s$. Combining~\eqref{102} with inequality~\eqref{103}
in which $z_0=U(0,s)z(0)$ and $t=s+1$, we obtain
\begin{align*}
|z(s+1)|_V^2&\leq C_3\biggl(|U(0,s)z(0)|_H^2+\int_{(s,s+1)}|f(\tau)|_H^2{\rm d}\tau\biggr)\\
&\le C_4e^{-\lambda
s}\biggl(|z(0)|_H^2+\sup_{t\ge0}\int_{(t,t+1)}e^{2\lambda
\tau}|f(\tau)|_H^2{\rm d}\tau\biggr).
\end{align*}
Using now~\eqref{104}, for $s\ge1$ we derive
\begin{equation} \label{106}
|z(s)|_V^2+\int_{(s,s+1)}|z(\tau)|_{{\rm D}(L)}^2{\rm d}\tau \le
C_5e^{-\lambda
s}\biggl(|z(0)|_H^2+\sup_{t\ge0}\int_{(t,t+1)}e^{2\lambda
\tau}|f(\tau)|_H^2{\rm d}\tau\biggr).
\end{equation}
On the other hand, it follows from~\eqref{104} that
\begin{equation} \label{107}
 \sup_{0\le s\le 1}|z(s)|_V^2+\int_{(0,1)}|z(\tau)|_{{\rm D}(L)}^2{\rm d}\tau
\le C_3\biggl(|z_0|_V^2+\int_{(0,1)}|f(\tau)|_H^2{\rm d}\tau\biggr).
\end{equation}
The required inequality~\eqref{est-b} follows immediately from~\eqref{106} and~\eqref{107}.

\medskip
{\it Step 2}.
We now prove that $\Xi$ maps the set $\mathcal Z^\lambda_\vartheta$ into itself.
Inequality~\eqref{est-b} with $f(t)=-Ba(t)$ implies that
\begin{equation} \label{108}
 |\Xi(a)|_{\mathcal Z^\lambda}^2\le C_1\biggl(|v_0|_V^2+\sup_{t\ge0}\int_{(t,t+1)}e^{2\lambda s}|Ba(s)|_H^2{\rm d}s\biggr).
\end{equation}
Now note that $|Ba|_H\le C_6|a|_V|a|_{{\rm D}(L)}$, whence it follows that
$$
\sup_{t\ge0}\int_{(t,t+1)}e^{2\lambda s}|Ba(s)|_H^2ds \le
C_6^2\sup_{t\ge0}\int_{(t,t+1)}\bigr(e^{\lambda
s}|a|_V^2\bigl)\,\bigr(e^{\lambda s}|a|_{{\rm D}(L)}^2\bigl){\rm
d}s \le C_6^2|a|_{\mathcal Z^\lambda}^4.
$$
Substituting this into~\eqref{108}, we see that if $a\in \mathcal Z^\lambda_\vartheta$, then
\begin{equation} \label{109}
|\Xi(a)|_{\mathcal Z^\lambda}\le C_7\bigl(|v_0|_V+|a|_{\mathcal Z^\lambda}^2\bigr)\le C_7\bigl(1+\vartheta |v_0|_V\bigr)\,|v_0|_V.
\end{equation}
Setting~$\vartheta=2C_7$ and choosing~$\epsilon>0$ so small that $C_7(1+\vartheta\epsilon)\le\vartheta$, we see that if $|v_0|_V\le\epsilon$, then~$\Xi$ maps the set~$\mathcal Z^\lambda_\vartheta$ into itself.

\medskip
{\it Step 3}.
It remains to prove that~$\Xi$ satisfies inequality~\eqref{contraction}.
Let us take two functions $a_1, a_2\in \mathcal Z^\lambda_\vartheta$ and set $a=a_1-a_2$ and $z=\Xi(a_1)- \Xi(a_2)$. Then
the function~$z$ satisfies the initial condition $z(0)=0$ and equation~\eqref{1000} with $f=Ba_2-Ba_1$.
Therefore, by inequality~\eqref{est-b}, we have
\begin{equation} \label{110}
|\Xi(a_1)-\Xi(a_2)|_{\mathcal Z^\lambda}^2\le
\sup_{t\ge0}\int_{(t,t+1)}e^{2\lambda s}|Ba_1-Ba_2|_H^2{\rm d}s.
\end{equation}
Using a standard estimate for $B(u,v)$ and the inequality $|u|_{L^\infty}^2\le C|u|_V|u|_{{\rm D}(L)}$, we derive
\begin{align*}
|Ba_1-Ba_2|_H^2&=|B(a_1,a)-B(a,a_2)|_H^2\\
&\leq C_8\bigl(|a_1|_{L^\infty}|a|_V+|a|_{L^\infty}|a_2|_V\bigr)^2\\
&\leq C_9\bigl(|a_1|_{V}|a_1|_{{\rm D}(L)}|a|_V^2+|a|_V|a|_{{\rm D}(L)}|a_2|_V^{2}\bigr).
\end{align*}
It follows that
\begin{equation} \label{111}
 \int_{(t,t+1)}e^{2\lambda s}|Ba_1-Ba_2|_H^2{\rm d}s
\le C_{10}\bigl(|a_1|_{\mathcal Z^\lambda}^2+|a_2|_{\mathcal Z^\lambda}^2\bigr) |a|_{\mathcal Z^\lambda}^2.
\end{equation}
Substituting~\eqref{111} into~\eqref{110} and recalling the definition of~$\mathcal Z^\lambda_\vartheta$, we obtain
$$
|\Xi(a_1)-\Xi(a_2)|_{\mathcal Z^\lambda}^2\le 2\vartheta C_{10}|v_0|_V^2|a_1-a_2|_{\mathcal Z^\lambda}^2.
$$
Choosing $\epsilon>0$ so small that $2\vartheta C_{10}\epsilon^2\le\gamma^2$, we see that if $|v_0|_V\le\epsilon$, then~\eqref{contraction} holds.
This completes the proof of the proposition.

\section{Appendix}
\label{S:appendix} 
\subsection{Karush--Kuhn--Tucker theorem}
Let~$\XX$ and~$\YY$ be Banach spaces and let $J:\XX\to\R$ and
$F:\XX\to\YY$ be two continuously differentiable functions.
Consider the following minimization problem with constraints:
\begin{equation} \label{Jmin-F}
J(x)\to \min, \quad F(x)=0.
\end{equation}
We will say that $\bar x\in\XX$ is a {\it local minimum\/}
for~\eqref{Jmin-F} if $F(\bar x)=0$ and there is a neighborhood
$U\ni \bar x$ such that $J(\bar x)\le J(x)$ for any $x\in U$ such that $F(x)=0$. A proof
of the following theorem can be found in~\cite{IT1979}.

\begin{theorem} \label{T:Jmin-F}
Let $\bar x\in \XX$ be a local minimum for~\eqref{Jmin-F} and let
the derivative $F'(\bar x):\XX\to\YY$ be a surjective operator.
Then there is $y^*\in\YY^*$ such that
\begin{equation} \label{kkt-cond}
J'(\bar x)+y^*\circ F'(\bar x)=0.
\end{equation}
\end{theorem}

\subsection{Quadratic functionals with linear constraint}
\label{quad-lin} Let~$\XX$ and~$\YY$ be normed vector spaces, let
$\tilde J(x,y)$ be a bounded symmetric bilinear form on~$\XX$ that is
weakly continuous with respect to each of its arguments, and let
$A:\XX\to\YY$ be a continuous surjective linear operator. Given a vector
$y\in\YY$, consider the minimization problem
\begin{equation} \label{Jmin-A}
J(x)\to \min, \quad Ax=y,
\end{equation}
where $J(x)=\tilde J(x,x)$. We will say that $\bar x\in\XX$ is a {\it
global minimum\/} for~\eqref{Jmin-A} if $A\bar x=y$ and $J(\bar
x)\le J(x)$ for $x\in \XX$ such that $Ax=y$. The following result is rather
standard in the optimal control theory, even though we were not able to find in the literature the statement we need.

\begin{theorem} \label{T:Jmin-A}
Suppose that $J(x)$ is non-negative and vanishes only for $x=0$,
and the set $\{x\in\XX:J(x)\le c\}$ is weakly compact for any
$c>0$. Then problem~\eqref{Jmin-A} has a unique global minimum
$\bar x\in\XX$, and the function $L:\YY\to\XX$ taking~$y$ to~$\bar
x$ is linear.
\end{theorem}

\begin{proof}
{\it Existence}. Let $m\ge0$ be the infimum of~$J$ on~$A^{-1}(y)$
and let~$\{x_n\}\subset A^{-1}(y)$ be a sequence such that
$J(x_n)\to m$. Since the set $\{x\in\XX:J(x)\le m+1\}$ is weakly
compact, we can assume that~$\{x_n\}$ converges weakly to a
vector~$\bar x\in \XX$. Now note that
$$
0\le J(x_n-\bar x)=J(x_n)-2\tilde J(x_n,\bar x)+J(\bar x).
$$
Combining this with the weak continuity of~$J$, we see that
$$
J(\bar x)\le\liminf_{n\to\infty} J(x_n)=m.
$$
Thus, $\bar x$ is a global minimum for~\eqref{Jmin-A}.

\smallskip
{\it Uniqueness}. Since the only point of~$\XX$ at which~$J$
vanishes is $x=0$, a standard argument proves that~$J$ is strictly
convex, that is,
$$
J\bigl(\tfrac{x_1+x_2}{2}\bigr)\le
\tfrac12\bigl(J(x_1)+J(x_2)\bigr)\quad\mbox{for all
$x_1,x_2\in\XX$},
$$
and the equality holds if and only if $x_1=x_2$. This immediately implies that the global minimum is
unique.

\smallskip
{\it Linearity}. Let $y\in\mathcal Y$ and $z\in A^{-1}(0)$. For
all $\lambda>0$, we have $A(Ly\pm\lambda z)=y$, and the definition
of~$L$ implies that $J(Ly)\leq {\tilde J}(Ly\pm\lambda z)$. It follows that
$0\leq\lambda J(z)\pm 2 {\tilde J}(Ly,\,z)$ for all $\lambda>0$. Letting
$\lambda$ go to $0$, we see that
\begin{equation} \label{100}
{\tilde J}(Ly,\,z)=0\quad\text{for all}\quad y\in\mathcal Y,\,z\in A^{-1}(0).
\end{equation}
For $a,\,b\in\mathcal Y$ and $\alpha,\,\beta\in\mathbb R$, let us set
$$
k:=\alpha La+\beta Lb -L(\alpha a+\beta b).
$$
Then $Ak=0$, and by~\eqref{100}, we have
$J(k)=\tilde J(k,\,k) 
=\alpha {\tilde J}(La,\,k)+\beta {\tilde J}(Lb,\,k) -{\tilde J}(L(\alpha a+\beta b),\,k)=0$.
It follows that $k=0$, and therefore~$L$ is linear.
\end{proof}

\subsection{Truncated observability inequality}
\label{s5.3} We first recall a well-known observability inequality
for the linearized Navier--Stokes system. Let us fix a function
$\hat u \in L^2(I_\tau,\,V)\cap\WW_\tau$, where~$\WW_\tau$ stands for the space of measurable vector-functions on~$Q_\tau$ such that (cf.\eqref{Wspace})
$$
|u|_{\WW_\tau}:=\sum_{j,\alpha}\,
\esssup_{(t,x)\in Q_\tau}\,\bigl|\p_t^j\p_x^\alpha u(t,x)\bigr|<\infty,
$$
where the sum is taken over $j=0,1$ and $|\alpha|\le1$. Consider the problem
\begin{align}
q_t-Lq-\mathbb B^\ast(\hat u)q&=0,\quad t\in I_{\tau},
\label{eq.q-gen}\\
q(\tau+1)&=q_1,\label{q1-gen}
\end{align}
where $q_1\in H$. By Theorem~2.2 in~\cite{iman01} (see also~\cite{FGIP-2004}), for any open subset
$\omega\subset\Omega$ there is a constant $C_\omega$ such that
\begin{equation}\label{obs-in}
|q(\tau)|_H^2\leq
C_\omega\int_{I_{\tau}}|q|_{L^2(\omega,\R^3)}^2\,{\rm d}t,
\end{equation}
Since $\supp\chi\cap\Omega\ne\emptyset$, the domain
$\omega_\chi:=\{x\in\Omega\mid\,|\chi(x)|>\rho\}$ is nonempty for a
sufficiently small~$\rho>0$. It follows from~\eqref{obs-in} that
$$
|q(\tau)|_H^2\leq
C_{\omega_\chi}\int_{I_{\tau}}|q|_{L^2(\omega_\chi,\R^3)}^2\,{\rm
d}t\leq C_{\omega_\chi}\rho^{-2}\int_{I_{\tau}}|\chi
q|_{L^2}^2\,{\rm d}t.
$$
Thus, setting $D_\chi^\prime:=C_{\omega_\chi}\rho^{-2}$, for any
solution of system \eqref{eq.q-gen},~\eqref{q1-gen}, we have the
observability inequality
\begin{equation}\label{obs-in.chi}
|q(\tau)|_H^2\leq D_\chi^\prime\int_{I_{\tau}}|\chi
q(t)|_{L^2}^2\,{\rm d}t.
\end{equation}
The following proposition shows that if~$q_1$ belongs to a
finite-dimensional subspace of~$H$, then the function~$\chi q$ on
the right-hand side of~\eqref{obs-in.chi} can be replaced
by~$P_M(\chi q)$ with a sufficiently large~$M$.

\begin{proposition}\label{obs-in.trunc}
For any $N\ge1$ there is an integer $M_1=\overline C_{[|\hat u|_{\mathcal W_\tau},N]}\in\N_0$ such that any solution~$q$ for system~\eqref{eq.q-gen}, \eqref{q1-gen} with $q_1\in F_N=\Pi_NH$ satisfies the inequality
\begin{equation}\label{obs-in.M}
 |q(\tau)|_H^2 \leq D_\chi
\int_{I_{\tau}}|P_{M_1}(\chi q(t))|_{L^2}^2\,{\rm d}t
\end{equation}
for a suitable constant $D_\chi$ depending only on $\chi$.
\end{proposition}

To prove the proposition, we need the following lemma.

\begin{lemma}\label{L:qV<qH}
For any solution $q$ of system \eqref{eq.q-gen},~\eqref{q1-gen} with $q_1\in F_N$, we have
\begin{equation} \label{120}
\int_{I_{\tau}}|\chi q(t)|_{H^1(\Omega,\R^3)}^2\,{\rm d}t 
\leq C \int_{I_{\tau}}|\chi q(t)|_{L^2(\Omega,\R^3)}^2\,{\rm d}t,
\end{equation}
where the constant $C$ depends only on~$N$, $\Omega$, and~$|\hat u|_\WW$.
\end{lemma}

\begin{proof}
We argue by contradiction. Suppose there is a sequence $(q^n_1,\,\hat u^n)\in F_N\times (L^2(I_{\tau},\,V)\cap \mathcal W_\tau)$,
with $(|\hat u^n|_{\mathcal W_\tau})$ bounded, such that the solution~$q^n$ of the problem
\begin{align}
q^n_t-Lq^n-\mathbb B^\ast(\hat u^n)q^n&=0,\quad t\in I_{\tau},
\label{eq.qn}\\
q^n(\tau+1)&=q^n_1\label{qn1}
\end{align}
satisfies the inequality
\begin{equation}\label{H1L2.abs-q}
\int_{I_{\tau}}|\chi q^n|_{H^1}^2\,{\rm
d}t>n\int_{I_{\tau}}|\chi q^n|_{L^2}^2\,{\rm d}t.
\end{equation}
Since the equations are linear, there is no loss of generality in
assuming that $|q^n_1|=1$. The boundedness of $(|\hat
u^n|_{\mathcal W_\tau})$ implies that~$(\p_x^\alpha\hat u^n)$ and~$(\p_x^\alpha\hat u^n_t)$ are bounded in~$L^\infty(Q_\tau)$ for $|\alpha|\le1$. It follows from Lemma~\ref{bd.hatu-r0} that the sequences $(q^n)$ and $(q^n_t)$ are bounded in $L^2(I_\tau,\,{\rm D}(L))$ and $L^2(I_\tau,\,H)$, respectively.
Since the unit ball in a Hilbert space is weakly compact and the unit ball in $L^\infty(Q_\tau)$ is compact in the weak$^*$ topology, there is a subsequence of $(q_1^n,q^n,\hat u^n)$ (for which we preserve the same notation),
a unit vector $q_1^\infty\in F_N$, and functions $q^\infty\in W(I_\tau,{\rm D}(L),H)$ and $\hat u^\infty\in\WW_\tau$ such that
\begin{align*}
q^{n}_1\quad&\to\quad q^\infty_1\quad\text{in}\quad F_N,\\
q^{n}\quad&\to\quad q^\infty\quad\text{in}\quad L^2(I_{\tau},\,V),\\
\partial_t q^{n}\quad&\rightharpoonup\quad\!\! \partial_t q^\infty\quad\text{in}\quad L^2(I_{\tau},\,H),\\
\hat u^{n}\quad&\to\quad\hat u^\infty\quad\text{in}\quad L^2(I_{\tau},\,H),\\
\partial_t^j\partial_x^\alpha\hat u^{n}\quad&\rightharpoonup_\ast\quad \partial_t^j\partial_x^\alpha\hat u^{\infty}
\quad\text{in}\quad L^\infty(Q_\tau),
\end{align*}
where $j=0,1$ and $|\alpha|\leq 1$.
Combining this with the boundedness of the sequences~$(\hat u^{n})$ and~$(q^{n})$ in the corresponding spaces,
we can easily pass to the limit in~\eqref{eq.qn}, \eqref{qn1} and derive the equations
\begin{align}
q^\infty_t-Lq^\infty-\mathbb B^\ast(\hat u^\infty)q^\infty&=0,\quad t\in I_{\tau},
\label{eq.pinf}\\
q^\infty(\tau+1)&=q^\infty_1.\label{pinf1}
\end{align}
Furthermore, since multiplication by~$\chi$ is a continuous operator in $L^2(I_\tau,\,H^1)$, we also have
\begin{equation}\label{lim.chi.p}
\chi q^n\to \chi q^\infty \quad\text{in}\quad L^2(I_{\tau},\,H^1(\Omega,\R^3)).
\end{equation}
Therefore, passing to the limit in inequality~\eqref{H1L2.abs-q} as $n\to\infty$, we conclude that
\begin{equation} \label{zero}
\int_{I_{\tau}}|\chi q^\infty|_{L^2}^2\,{\rm d}t=0.
\end{equation}
Applying now the observability inequality~\eqref{obs-in.chi} to equation~\eqref{eq.pinf}
considered on the interval $(\tau+r,\tau+1)$ with $0\le r<1$,
we conclude that $q^\infty(t)=0$ for $\tau \le t<\tau+1$. Since $q^\infty\in C(\bar I_\tau,V)$, we obtain $q_1^\infty=q^\infty(\tau+1)=0$. This contradicts the fact that $q_1^\infty\in F_N$ is a unit vector. The contradiction obtained proves that~\eqref{120} holds.
\end{proof}

{\em Proof of Proposition \ref{obs-in.trunc}}. 
We use Lemma~\ref{L:qV<qH}  to derive
\begin{align*}
\int_{I_{\tau}}|\chi q|_{L^2}^2\,{\rm d}t
&\leq \int_{I_{\tau}}|P_M(\chi q)|_{L^2}^2\,{\rm d}t+\int_{I_{\tau}}|(1-P_M)\chi q|_{L^2}^2\,{\rm d}t\\
&\leq \int_{I_{\tau}}|P_M(\chi q)|_{L^2}^2\,{\rm d}t+\beta_M^{-1}\int_{I_{\tau}}|(1-P_M)(\chi q)|_{H^1}^2\,{\rm d}t\\
&\leq \int_{I_{\tau}}|P_M(\chi q)|_{L^2}^2\,{\rm d}t+\beta_M^{-1}\int_{I_{\tau}}|\chi q|_{H^1}^2\,{\rm d}t\\
&\leq \int_{I_{\tau}}|P_M(\chi q)|_{L^2}^2\,{\rm
d}t+\beta_M^{-1}\overline C_{[N,|\hat u|_{\mathcal
W_\tau}]}\int_{I_{\tau}}|\chi q|_{L^2}^2\,{\rm d}t.
\end{align*}
Recall that $\beta_j$ stands for the $j^{\text{th}}$ eigenvalue of the Dirichlet Laplacian.
Choosing the integer~$M=M_1$ so large that $\beta_{M_1}^{-1}\overline C_{[N,|\hat u|_{\mathcal W_\tau}]}\leq\frac{1}{2}$, we obtain
\begin{equation*}
\int_{I_{\tau}}|\chi q|_{L^2}^2\,{\rm d}t\leq
2\int_{I_{\tau}}|P_{M_1}(\chi q)|_{L^2}^2\,{\rm d}t.
\end{equation*}
Combining this with~\eqref{obs-in.chi}, we arrive the required inequality~\eqref{obs-in.M}.
\endproof

\medskip
\thanks{{\bf Acknowledgments.} This work was supported by LEA CNRS Franco-Roumain ``Math\'ematiques \& Mod\'elisation'' 
in the framework of the project {\it Control of nonlinear PDE's}. We thank the anonymous referees for pertinent critical remarks that helped to improve the presentation and to eliminate some inaccuracies of the previous version of the paper.}
 

\providecommand{\bysame}{\leavevmode\hbox to3em{\hrulefill}\thinspace}
\providecommand{\MR}{\relax\ifhmode\unskip\space\fi MR }
\providecommand{\MRhref}[2]{%
  \href{http://www.ams.org/mathscinet-getitem?mr=#1}{#2}
}
\providecommand{\href}[2]{#2}

\end{document}